\documentclass[11pt]{article}
\usepackage{graphicx}

\def\beq{\begin{eqnarray}}
\def\eeq{\end{eqnarray}}

\begin{document}

\fontsize{11}{14.5pt}\selectfont

\vspace*{1in}

\begin{center}

\begin{center} \Large \bf 
  Fast Exact Summation Using Small and Large Superaccumulators
\end{center}

\vspace{8pt}

{\large Radford M. Neal \\[4pt]
  \normalsize Dept.\ of Statistical Sciences and Dept.\ of Computer Science \\
  University of Toronto \\[4pt]
  \texttt{http://www.cs.utoronto.ca/$\sim$radford/}\\
  \texttt{radford@stat.utoronto.ca}\\[4pt]
  20 May 2015}
 
\end{center}

\vspace{20pt}

\noindent I present two new methods for exactly summing a set of
floating-point numbers, and then correctly rounding to the nearest
floating-point number. Higher accuracy than simple summation (rounding
after each addition) is important in many applications, such as
finding the sample mean of data. Exact summation also guarantees
identical results with parallel and serial implementations, since the
exact sum is independent of order. The new methods use variations on
the concept of a ``superaccumulator'' --- a large fixed-point number
that can exactly represent the sum of any reasonable number of
floating-point values.  One method uses a ``small'' superaccumulator
with sixty-seven 64-bit chunks, each with 32-bit overlap with the next
chunk, allowing carry propagation to be done infrequently. The small
superaccumulator is used alone when summing a small number of
terms. For big summations, a ``large'' superaccumulator is used as
well. It consists of 4096 64-bit chunks, one for every possible
combination of exponent bits and sign bit, plus counts of when each
chunk needs to be transferred to the small superaccumulator. To add a
term to the large superaccumulator, only a single chunk and its
associated count need to be updated, which takes very few instructions
if carefully implemented. On modern 64-bit processors, exactly summing
a large array using this combination of large and small
superaccumulators takes less than twice the time of simple, inexact,
ordered summation, with a serial implementation.  A parallel
implementation using a small number of processor cores can be expected
to perform exact summation of large arrays at a speed that reaches the
limit imposed by memory bandwidth. Some common methods that attempt to
improve accuracy without being exact may therefore be pointless, at
least for large summations, since they are slower than computing the
sum exactly.

\newpage

\subsection*{Introduction}\vspace*{-7pt}

\noindent Computing the sum of a set of numbers can produce an
inaccurate result if it is done by adding each number in turn to an
accumulator with limited precision, with rounding performed on each
addition.  The final result can be much less accurate than the
precision of the accumulator if cancellation occurs between positive
and negative terms, or if accuracy is lost when many small numbers are
added to a larger number. Such inaccuracies are a problem for many
applications, one such being the computation of the sample mean of
data in statistical applications.

Much work has been done on trying to improve the accuracy of
summation.  Some methods aim to somewhat improve accuracy at little
computational cost, but do not guarantee that the result is the
correctly rounded exact sum.  For example, Kahan's method (Kahan,
1965) tries to compensate for the error in each addition by
subtracting this error from the next term before it is added.  Another
simple method is used by the R language for statistical computation
(R~Core Team, 1995--2015), which computes the sample mean of data by first
computing a tentative mean (adding terms in the obvious way) and then
adjusting this tentative mean by adding to it the mean (again,
computed in the obvious way) of the difference of each term from the
tentative mean. This method sometimes improves accuracy, but can also
make the result less accurate. (For example, the R expression
\verb|mean(c(1e15,-1e15,0.1))| gives a result accurate to only three
decimal digits, whereas the obvious method would give the exact mean
rounded to about 16 digits of accuracy).

Many methods have been developed that instead compute the
\textit{exact} sum of a set of floating-point values, and then
correctly round this exact sum to the closest floating-point value.
This obviously would be preferable to any non-exact method, if the
exact computation could be done sufficiently quickly.

An additional advantage of exact methods is that they can easily be
parallelized, without changing the result, since unlike inexact
summation, the exact sum does not depend on the order in which terms
are added.  In contrast, parallelizing simple summation in the obvious
way, by splitting the sum into parts that are summed (inexactly) in
parallel, then adding these partial sums, will in general produce a
different result than the simple serial method.  Furthermore, the
result obtained will depend on the details of the parallel algorithm,
and perhaps on the run-time availability of processor cores.

Differing results will also arise from serial implementations that do
not sum terms in the usual left-to-right order.  Such implementations
are otherwise attractive, since many modern processors have multiple
computational units that can be exploited via instruction-level
parallelism, if data dependencies allow it.  Summing four numbers as
$((a_1+a_2)+a_3)+a_4$ does not allow for any parallelism, but summing
them as $(a_1+a_2)+(a_3+a_4)$ does, although it may produce a
different result.  In contrast, focusing on exact computation ensures
that any improvements in computational methods will not lead to
non-reproducible results.

Exact summation methods fall into two classes --- those implemented
using standard floating-point arithmetic operations available in
hardware on most current processors, such as the methods of Zhu and
Hayes (2010), and those that instead perform the summation with
integer arithmetic, using a ``superaccumulator''.  Hybrid methods
using both techniques been investigated by Collange, Defour, Graillat,
and Iakymchuk (2015a,b).

The methods of this paper can be seen as using ``small'' and ``large''
variations on a superaccumulator, though the ``large'' variation
resembles other superaccumulator schemes only distantly.  The general
concept of a superaccumulator is that it is a fixed-point numerical
representation with enough binary digits before and after the binary
point that it can represent the sum of any reasonable number of
floating-point values exactly and without overflow.  Such a scheme is
possible because the exponent range in floating number formats is
limited.

The idea of such a superaccumulator goes back at least to Kulisch and
Miranker (1984), who proposed its use for exact computation of dot
products.  In that context, the superaccumulator must accommodate the
range of possible exponents in a product of two floating-point
numbers, which is twice the exponent range of a single floating-point
number, and the terms added to the superaccumulator will have twice
the precision of a single floating-point number.  In this paper, I
will consider only the problem of summing individual floating-point
values, in the standard (IEEE Computer Society, 2008) 64-bit ``double
precision'' floating-point format, not higher-precision products of
such values.  Directly extending the methods in this paper to such
higher precision sums would require doing arithmetic with 128-bit
floating-point and integer numbers, which at present is typically
unsupported or slow.  In the other direction, exact dot products of
``single-precision'' (32-bit) floating-point values could be computed
using the present implementation, and exact summation of
single-precision values could be done even more easily (with smaller
superaccumulators).

Below, I first describe the standard floating-point and integer numeric
formats assumed by the methods of this paper, and then present the
``small'' superaccumulator method, whose design incorporates a
tradeoff between the largely fixed time for initialization and
termination and the additional time used for every term added.  I then
present a method in which such a small superaccumulator is combined
with a ``large'' superaccumulator.  This method has a higher fixed
cost, but requires less time per term added.

I evaluate the performance of the small and large superaccumulator
methods using a carefully written implementation in C, which is
provided as supplementary information to this paper. \mbox{I~compare}
the performance of these new methods with the obvious (inexact) simple
summation method, with a variation on simple summation that
accumulates sums of of terms with even and odd indexes separately,
allowing for increased instruction-level parallelism, and with the
exact iFastSum and OnlineExact methods of Zhu and Hayes (2010), who
have provided a C++ implementation.  Timing tests are done on sixteen
computer systems, that use Intel, AMD, ARM, and Sun processors
launched between 2000 and 2013.

The results show that on modern 64-bit processors, when summing many
terms (tens of thousands or more), the large superaccumulator method
is less than a factor of two slower than simple inexact summation, and
is significantly faster than all the other exact methods tested.  When
summing fewer than about a thousand terms, the small superaccumulator
method is faster than the large superaccumulator method.  The iFastSum
method is almost always slower than the small superaccumulator method,
except for very small summations (less than about twenty terms), for
which it is sometimes slightly faster.  The OnlineExact method is
about a factor of two slower than the large superaccumulator method on
modern 64-bit processors.  It is also slower or no faster on older
processors, with the exception of 32-bit processors based on the
Pentium 4 Netburst architecture, for which it is about a factor of two
faster than the large superaccumulator method.

I conclude by discussing the implications of these performance
results, and the possibility for further improvements, such as methods
designed for small summations (less than 100 terms), methods using
multiple processor cores, and implementations of methods in
carefully-tuned assembly language.

\subsection*{Floating-point and integer formats}\vspace*{-7pt}

The methods in this paper are designed to work with floating-point
numbers in the standard (IEEE Computer Society, 2008) 64-bit
``double-precision'' format, which is today universally available, in
hardware implementations, on general-purpose computers, and used for
the C language \verb|double| data type.

\begin{figure}

\centerline{~\includegraphics[scale=0.5]{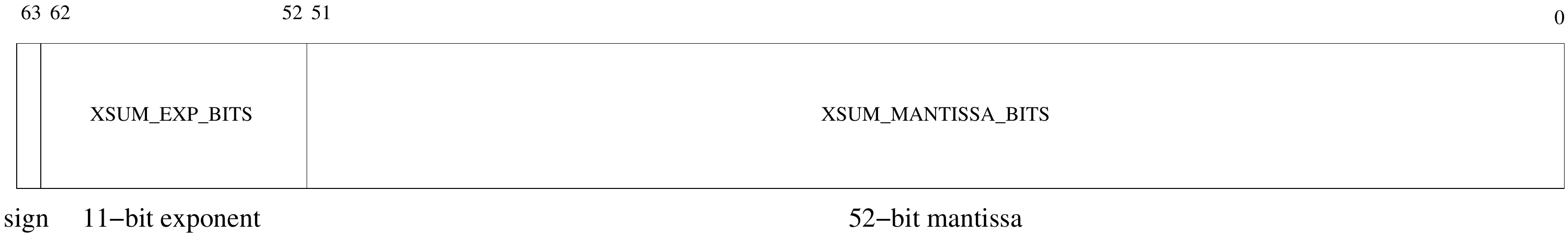}}

\vspace{8pt}

{\small

Note:\ \ The C code shown later uses the symbols \verb|XSUM_EXP_BITS|
(11) and \verb|XSUM_MANTISSA_BITS| (52), as well as the symbols
\verb|XSUM_EXP_MASK|, equal to \verb|(1<<XSUM_EXP_BITS)-1|, and
\verb|XSUM_MANTISSA_MASK|, equal to
\verb|((int64_t)1<<XSUM_MANTISSA_BITS)-1|.

}

\caption{Format of an IEEE 64-bit floating-point number.}\label{fig-ieee}

\end{figure}

Numbers in this format, illustrated in Figure~\ref{fig-ieee}, consist
of a sign bit, $s$, 11 exponent bits, $e$, and 52 mantissa bits, $m$.
Interpreting each group of bits as an integer in binary notation, when
$e$ is not 0 and not 2047, the number represented by these bits is
$(-1)^s \,\times\, 2^{e-1023} \,\times\, (1\!+\!m2^{-52})$. That
is, $e$ represents the binary exponent, with a bias of 1023, and the
full mantissa consists of an implicit~1 followed by the bits of $m$.
When $e$ is 0 (indicating a ``denormalized'' number), the number
represented is $(-1)^s \,\times\, 2^{-1022} \,\times\, m2^{-52}$.
That is, the true exponent is $1\!-\!1023$, and the mantissa does not
include an implicit~1.  A value for $e$ of 2047 indicates plus or
minus ``infinity'' when $m$ is zero, and a special NaN (``Not a
Number'') value otherwise.  Note that the smallest non-zero
floating-point value is $2^{-1074}$ and the largest non-infinite value
is $2^{1023}\,\times\,(2-2^{-52})$.

I also assume that unsigned and signed (two's complement) 64-bit
integer formats are available, and are accessible from C by the
\verb|uint64_t| and \verb|int64_t| data types.  These formats are
today universally available for general purpose computers, and
accessible from C in implementations compliant with the C99 standard.
Arithmetic on 64-bit quantities is well-supported by recent 64-bit
processors, but even on older \mbox{32-bit} processors, addition,
subtraction, and shifting of 64-bit quantities are not extraordinarily
slow, being facilitated by instructions such as ``add with carry''.

Finally, I assume that the byte ordering of 64-bit floating-point
values and 64-bit integers is consistent, so that a C union type with
\verb|double|, \verb|int64_t|, and \verb|uint64_t| fields will allow
access to the sign, exponent, and mantissa of a 64-bit floating-point
value stored into the \verb|double| field via shift and mask
operations on the 64-bit signed and unsigned integer fields.  Such
consistent ``endianness'' is not guaranteed by any standard, but
seems to be nearly universal on today's computers (of both ``big
endian'' and ``little endian'' varieties) --- including Intel x86,
SPARC, and modern ARM processors (though it appears some past ARM
architectures may not have been consistent).

\subsection*{Exact summation using a small superaccumulator}\vspace*{-7pt}

I first present a new summation method using a relatively small
superaccumulator, which will prove to be the preferred method for
summing a moderate number of terms, and which is also a component of
the large superaccumulator method presented below.  The details of
this scheme are designed for fast implementation in software, in
contrast to some other designs (eg, Kulisch, 2011) that are meant
primarily for hardware implementation.

The most obvious design of a superaccumulator for use in summing
64-bit floating-point values would be a fixed-point binary number
consisting of a sign bit, $1024+\lceil\log_2 N\rceil$ bits to the
left of the binary point, where $N$ is the maximum number of terms
that might be summed, and $1074$ bits to the right of the binary
point.  The bits of such a superaccumulator could be stored in around
34 64-bit words.

However, this representation has several disadvantages.  When adding a
term to the superaccumulator, carries might propagate through several
64-bit words, requiring a loop in the time-critical addition
operation.  Furthermore, this sign-magnitude representation requires
that addition and subtraction be handled separately, with the sign
changing as necessary, necessitating additional complexities.  If the
superaccumulator instead represents negative numbers in two's
complement form, additions that change the sign of the sum will need
to alter all the higher-order bits.

Carry propagation can be sped up using a somewhat redundant
``carry-save'' representation, in which the high-order bits of each
64-bit ``chunk'' of the superaccumulator overlap the low-order bits of the
next higher chunk, allowing carry propagation to be deferred for some
time.  This approach is used, for example, by Collange, et al.\
(2015a,b), whose chunks have 8-bit overlap.  In the scheme of
Collange, et al., chunks can apparently also have different signs, an
arrangement that can alleviate the problems of representing negative
numbers, by allowing local updates without the need to determine the
overall sign of the number immediately.

\begin{figure}

\centerline{~\includegraphics[scale=0.5]{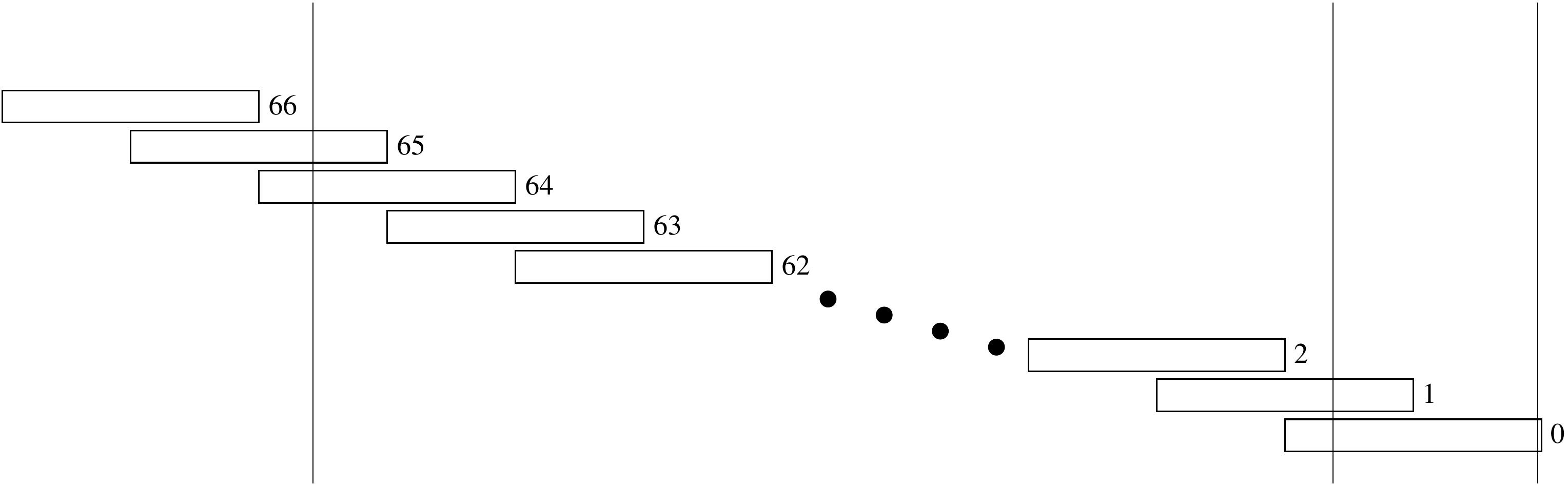}}

\vspace{8pt}

{\small

Note:\ \ The C code shown later uses the symbols \verb|XSUM_SCHUNKS| (67),
\verb|XSUM_LOW_MANTISSA_BITS| (32), \verb|XSUM_HIGH_EXP_BITS| (6) and
\verb|XSUM_LOW_EXP_BITS| (5), along with corresponding masks.

}

\caption{Chunks making up a small superaccumulator. There are 67
chunks in the superaccumulator, whose indexes (shown to the right) are
related to the high 6 bits of the exponent in a number, with the low 5
bits of an exponent specifying a position within a chunk.  Each chunk
is a 64-bit signed integer, with chunks overlapping by 32 bits.
Chunks are shown with overlap above, so that horizontal position
corresponds to the positional value of each bit. The vertical lines at
the right delimit the range of denormalized numbers (note that the
rightmost bit is unused).  The vertical line at the left is the
position of the topmost implicit 1 bit of the largest possible 64-bit
floating point number.  Bits to the left of that are provided to
accomodate larger numbers that can arise when many numbers are
summed.}\label{fig-small-acc}

\end{figure}

In the design I use here, the small superaccumulator consists of 67 
signed (two's complement) 64-bit chunks, with 32-bit overlap.  Chunks
are indexed starting at 0 for the lowest-order chunk.  
Denoting the value of chunk $i$ as $c_i$,
the number represented by the superaccumulator is defined to be
\[
   \sum_{i=0}^{66} \,c_i\, 2^{32i-1075}
\]
The $c_i$ will always be in the range $-(2^{63}\!-\!1)$ to
$2^{63}\!-\!1$. For convenience, the representation is further
restricted so that the highest-order chunk (for $i=66$) is in the
range $-2^{32}$ to $2^{32}\!-\!1$.  This representation is diagrammed
in Figure~\ref{fig-small-acc}.

The largest number representable in this superaccumulator is
$2^{1069}-2^{-1074}$.  It can therefore represent any sum of up to
$2^{45}$ terms, which would occupy more than 281 terabytes of memory.
This capacity to represent values beyond the exponent range of the
64-bit floating-point format ensures that the final rounded 64-bit
floating-point sum obtained using the superaccumulator will be finite
whenever the final exact sum is within range, even when summing the
values in the ordinary way would have produced overflow for an
intermediate result.  This is an advantage over methods such as those
of Zhu and Hayes (2010), which use floating-point arithmetic, and
hence cannot bypass temporary overflows.

Due to the overlap of chunks, and the possibility that they have
different signs, a single number can have many possible
representations in the superaccumulator.  However, a canonical form
is produced when carry propagation is done, which happens periodically
when adding terms to the superaccumulator, and whenever a
floating-point number that is the correct rounding of the
superaccumulator's value is needed.  Carry propagation starts at the
low order chunk ($i=0$), and proceeds by clearing the high-order
32-bits of each chunk to zero, and adding these bits (regarded as a
signed integer) to the next-higher chunk.  The process ends when we
reach the highest-order chunk, whose high-order 32 bits will be either
all 0s or all 1s, depending on whether the number is positive or
negative.  Note that all chunks other than this highest-order chunk
are positive after carry propagation.

If carried out as just described, carry propagation for a negative
number could require modification of many higher-order chunks, all of
which would be set to $-1$ (ie, all 1s in two's complement).  To avoid
this inefficiency, the procedure is modified so that such high-order
chunks that would have value $-1$ are instead set to zero, and the
upper 32-bits of the next-lower chunk are set to all 1s (so that it is
now negative), which produces the same represented number.

After carry propagation, all chunks will be no larger than $2^{32}$ in
absolute value.  In the procedure described next for adding a
floating-point value to the superaccumulator, the amount added to (or
subtracted from) any chunk is at most $2^{52}\!-\!1$.  It follows that
the values of all chunks are guaranteed to remain within their allowed
range if no more than $2^{11}\!-\!1 = 2047$ additions are done between
calls of the carry propagation routine.  This is sufficiently large
that it makes sense to keep only a global count of remaining additions
before carry propagation is needed, rather than keeping counts for
each chunk, or detecting actual overflow when adding to or subtracting
from a chunk.  Using only a global count will result in carry
propagation being done more often than necessary, but since the cost
of carry propagation should be only a few tens of instructions per
chunk, reducing calls to the carry propagation routine cannot justify
even one additional instruction in the time-critical addition
procedure.

Addition of a 64-bit floating-point value to the superaccumulator
starts with extraction of the 11 exponent bits and 52 mantissa bits,
using shift and mask operations that treat the value as a 64-bit
integer. Note that the sign of the floating-point number is the same
as the sign of its 64-bit integer form, so no extraction of the sign
bit is necessary.

If the exponent bits are all 1s, the floating-point value is an
infinity or a NaN, which are handled specially by storing indicators in
auxiliary Inf and NaN fields of the superaccumulator.  This operation
is typically not highly time-critical, since Inf and NaN operands are
expected to be fairly infrequent.

If the exponent bits are all 0s, the floating-point value is a zero or
a non-zero denormalized number.  If it is zero, the addition operation
is complete, since nothing need be done to add zero.  Otherwise, the
exponent is changed to 1, since this is the true exponent (with bias)
of denormalized numbers.

If the exponent bits are neither all 0s nor all 1s, the value is an
ordinary normalized number.  In this case, the implicit~1 bit that is
part of the mantissa value is explicitly set, so that the mantissa
value now contains 53 bits.

Further shift and mask operations separate the exponent into its
high-order 6 bits and low-order 5 bits.  The high-order exponent bits,
denoted $i$, index one of the first 64 chunks of the superaccumulator.
Chunks $i$ and $i+1$ will be modified by adding or subtracting bits of
the mantissa.  Due to the overlap of these chunks, this could be done
in several ways, but it seems easiest to modify chunk $i$ by adding or
subtracting a 32-bit value, and to use the remaining bits to modify
chunk $i+1$.

In detail, the quantity to add to or subtract from chunk $i$ is found
by shifting the 53-bit mantissa left by the number of bits given by
the low-order 5 bits of the exponent, and then masking out only the
low-order 32 bits.  The shift positions these mantissa bits to their
proper place in the superaccumulator.  The quantity to add to or
subtract from chunk $i+1$ is found by shifting the 53-bit mantissa
right by 32 minus the amount of the previous shift.  This isolates
(without need of a masking operation) the bits that were not used to
modify chunk $i$, positioning them properly for adding to or
subtracting from chunk $i+1$.  Note that this quantity will have at
most 52 bits, since at least 1 mantissa bit will be used to modify
chunk $i$.

When modifying both chunk $i$ and chunk $i+1$, whether to add or
subtract is determined by the sign of the number being added.  Note
that it is quite possible for different chunks to end up with
different signs after several terms have been added, but the overall
sign of the number is resolved when carry propagation is done.

The C code used for this addition operation is shown in
Figure~\ref{fig-small-add}.  A function that sums an array would use
this code (expanded from an inline function) in its inner loop that
steps through array elements.  This summation function must call the
carry propagation routine after every 2047 additions.  This is most
easily done with nested loops, with the inner loop adding numbers
until some limit is reached, which is the same form as the inner loop
would be if no check for carry propagation were needed.  

\begin{figure}

\begin{verbatim}
/*** Declarations of types used to define the small superaccumulator ***/

typedef int64_t xsum_schunk;     /* Integer type of small accumulator chunk */

typedef struct                   /* A small superaccumulator */
{ xsum_schunk chunk[XSUM_SCHUNKS]; /* Chunks making up small accumulator */
  int64_t Inf;                     /* If non-zero, +Inf, -Inf, or NaN */
  int64_t NaN;                     /* If non-zero, a NaN value with payload */
  int adds_until_propagate;        /* Number of remaining adds before carry */
} xsum_small_accumulator;          /*     propagation must be done again    */

/*** Code for adding the double 'value' to the small accumulator 'sacc' ***/

union { double fltv; int64_t intv; } u;

u.fltv = value;
ivalue = u.intv;
mantissa = ivalue & XSUM_MANTISSA_MASK;
exp = (ivalue >> XSUM_MANTISSA_BITS) & XSUM_EXP_MASK;

if (exp != 0 && exp != XSUM_EXP_MASK) /* normalized */
{ mantissa |= (int64_t)1 << XSUM_MANTISSA_BITS;
}
else if (exp == 0) /* zero or denormalized */
{ if (mantissa == 0) return;
  exp = 1;
}
else /* Inf or NaN */
{ xsum_small_add_inf_nan (sacc, ivalue);
  return;
}

low_exp = exp & XSUM_LOW_EXP_MASK;
high_exp = exp >> XSUM_LOW_EXP_BITS;
chunk_ptr = sacc->chunk + high_exp;
chunk0 = chunk_ptr[0];
chunk1 = chunk_ptr[1];
low_mantissa = (mantissa << low_exp) & XSUM_LOW_MANTISSA_MASK;
high_mantissa = mantissa >> (XSUM_LOW_MANTISSA_BITS - low_exp);

if (ivalue < 0)
{ chunk_ptr[0] = chunk0 - low_mantissa;
  chunk_ptr[1] = chunk1 - high_mantissa;
}
else
{ chunk_ptr[0] = chunk0 + low_mantissa;
  chunk_ptr[1] = chunk1 + high_mantissa;
}
\end{verbatim}

\caption{Extracts from C code for adding a 64-bit floating point value
         to a small superaccumulator.}\label{fig-small-add}

\end{figure}

Once all terms have been added to the small superaccumulator, a
correctly rounded value for the sum can be obtained, after first
performing carry propagation.  Special Inf and NaN values must be
handled specially.  Otherwise, the chunks are examined starting at the
highest-order non-zero chunk, and proceeding to lower-order chunks as
necessary.  Note that the sign of the rounded value is given by the
sign of the highest-order chunk.

Denormalized numbers are easy to identify, and do not require rounding.

For normalized numbers, a tentative exponent for the rounded value can
be obtained by converting the highest chunk's integer value to
floating point, and then looking at the exponent of the converted
value.  This is the only use of a floating-point operation in the
superaccumulator routines.  If desired, this operation could be
replaced with some other method of finding the topmost 1 bit in a
32-bit word (for instance, binary search using masks).  This tentative
exponent allows construction of a tentative mantissa from the
highest-order chunk and the next lower one or two chunks.  Chunks of
lower order may need to be examined in order to produce a correctly
rounded result, potentially all the way to the lowest-order chunk.
Rounding may change the final exponent.

See the code in the supplemental information for further (somewhat
finicky) details of rounding.  At present, only the commonly-used
``round to nearest, with ties to even'' rounding mode is implemented,
but implementing other rounding modes would be straightforward.

Exact summation using the small superaccumulator has a fixed cost, due
to the need to set all 67 chunks to zero initially, and to scan all
chunks when carry propagating in order to produce the final rounded
result.  As will be seen from the experiments below, this fixed cost
is roughly 12.5 times the cost of adding a single term to the
superaccumulator.  A naive count of operations in the C code of
Figure~\ref{fig-small-add} gives about 19 operations to add a term to
the superaccumulator, compared to 2 operations (fetch and add) for
simple floating-point summation. The actual per-term time ratio is not
that bad on modern 64-bit processors, probably because these
processors can exploit instruction-level parallelism. Nevertheless, to
obtain good performance for large summations, we are motivated to look
for a scheme with smaller cost per term, even if this increases fixed
overhead.

\subsection*{Faster exact summation of many terms with a large 
             superaccumulator}\vspace*{-7pt}

To reduce the per-term cost of summing values with a superaccumulator,
we would like to eliminate from the inner summation loop the
operations of testing for special Inf or NaN values, checking the sign
of the term in order to decide whether to add or subtract, and
splitting the mantissa bits into two parts, so they can be added to
different chunks. This can be accomplished by using a large
superaccumulator that has 4096 64-bit chunks, one for every possible
combination of sign and exponent bits, as well as 4096 16-bit counts,
one for each chunk.  We still use a small superaccumulator as well,
transferring partial sums from the large superaccumulator to the small
superaccumulator as necessary to avoid loss of information from
overflow.  The counts in the large superaccumulator are all initialized to
$-1$; the chunks are not set initially.

The C code for adding a value to this large superaccumulator is shown
in Figure~\ref{fig-large-add}.  It starts by isolating the sign and
exponent bits of the floating-point value, viewed as an unsigned
64-bit integer, by doing a right shift by 52 bits (with zero fill, so
no masking is needed).  These 12 bits will be used to index a 64-bit
chunk of the large superaccumulator, and the corresponding 16-bit
count.

\begin{figure}[t]

\begin{verbatim}
/*** Declarations of types used to define the large superaccumulator ***/

typedef uint64_t xsum_lchunk;      /* Integer type of large accumulator chunk,
                                      must be EXACTLY 64 bits in size */

typedef int_least16_t xsum_lcount; /* Signed int type of counts for large acc. */

typedef uint_fast64_t xsum_used;   /* Unsigned type for holding used flags */

typedef struct
{ xsum_lchunk chunk[XSUM_LCHUNKS]; /* Chunks making up large accumulator */
  xsum_lcount count[XSUM_LCHUNKS]; /* Counts of # adds remaining for chunks,
                                        or -1 if not used yet or special. */
  xsum_used chunks_used[XSUM_LCHUNKS/64];  /* Bits indicate chunks in use */
  xsum_used used_used;             /* Bits indicate chunk_used entries not 0 */
  xsum_small_accumulator sacc;     /* The small accumulator to condense into */
} xsum_large_accumulator;

/*** Code for adding the double 'value' to the large accumulator 'lacc' ***/

union { double fltv; uint64_t uintv; } u;

u.fltv = value
ix = u.uintv >> XSUM_MANTISSA_BITS;
count = lacc->count[ix] - 1;

if (count < 0)
{ xsum_large_add_value_inf_nan (lacc, ix, u.uintv);
}
else
{ lacc->count[ix] = count;
  lacc->chunk[ix] += u.uintv;
}
\end{verbatim}

\caption{Extracts from C code for adding a 64-bit floating point value
         to a large superaccumulator.}\label{fig-large-add}

\end{figure}

The count indexed by the sign and exponent is then fetched, and
decremented.  If this decremented count is non-negative, it is stored
as the new value for this count, and the entire floating-point value
is added, as a 64-bit integer, to the 64-bit chunk indexed by the sign
and exponent.  Note that no operation to mask out just the mantissa
bits of this value is done.  This masking can be omitted because the
undesired bits at the top are the same for every add to any particular
chunk, and are known from the index of that chunk. Since the number of
adds that have been done is also kept track of in the count, the
effect of adding these bits can be undone before transferring the sum
held in the chunk to the small superaccumulator.

If instead the decremented count is negative, a routine to do special
processing is called.  This test merges a check for the value being an
Inf or NaN, a check for the indexed chunk having not yet been
initialized, and a check for having already done the maximum allowed
number (4096) of adds to the indexed chunk, so that the chunk's
contents must now be transferred to the small superaccumulator.

Since all these circumstances are expected to arise infrequently, the
special processing routine is not time-critical.  It operates as
follows.  When the exponent bits in the index passed to this routine
are all 1s, the value being added is an Inf or NaN, which is handled
by setting special fields of the small superaccumulator associated
with this large superaccumulator.  The count for this index remains at
$-1$, so that subsequent adds of this Inf or NaN will also be
processed specially.  For other exponents, if the count (before being
decremented) is $-1$, indicating that this is the first use of this
chunk, the chunk is initialized to zero, and the count is set to 4096.
Otherwise, the count must be zero, indicating that the maximum of 4096
adds have previously been done to this chunk, in which case the sum is
transferred to the small superaccumulator, the chunk is reset to zero,
and the count is reset to 4096.  In the latter two cases, the addition
then proceeds as usual (adding to the chunk and decrementing the
count).

The partial sum in a large superaccumulator chunk will need to be
transferred to the small superaccumulator when the maximum number of
adds before overflow has already been done, or when the final rounded
result is desired.  When the maximum of 4096 adds has been done, the
bits in the chunk are the correct sum of mantissa bits, without any
further adjustment, since adding the same sign and exponent bits 4096
times is the same as multiplying by 4096, which is the same as
shifting these bits left 12 positions, which removes them from the
64-bit word.  When the transfer to the small superaccumulator is done
before 4096 adds to the chunk, we need to add to the chunk the the
chunk's index (the sign and exponent bits) times the count of
remaining allowed adds, shifted left 52 bits, which has the effect of
leaving only the sum of mantissa bits.

The sum of the mantissa bits for all values that were added to this
chunk has unsigned magnitude up to $2^{64}\!-\!2^{12}$, so all 64 bits
of the chunk are used.  There would be several ways of transferring
these bits to the small superaccumulator, but it seems easiest to do
so by modifying three consecutive small superaccumulator chunks by
adding or subtracting 32-bit quantities.  Conceptually, these three
32-bit quantities are obtained by shifting the 64-bit chunk left by
the number of positions given by the low 5 bits of the exponent (the
same as the low 5 bits of the chunk index), and then extracting the
lowest 32 bits, the next 32 bits, and the highest 32 bits.  However,
since shift operations on quantities greater than 64 bits in size may
not be available, the equivalent result is instead found using a some
left and some right shifts, and suitable masking operations.  For
chunks corresponding to normalized floating-point values (ie, for
which the exponent is not zero), we also add in the sum of all the
implicit~1 bits at the top of the mantissa (which would be beyond the
top of the 64-bit chunk) to the appropriate 32-bit quantity.  Finally
we either add or subtract these three 32-bit quantities from the
corresponding chunks of the small superaccumulator according to the
sign bit, which is the top bit of the 12-bit index of the chunk.

The fixed cost of summation using a large superaccumulator is greater
than that of using only a small superaccumulator because of the need
to initialize the array of 4096 counts, occupying 8192 bytes.  Note
this is in addition to the fixed costs of using the small
superaccumulator, which is still needed as well.  Note, though, that
the 4096 large superaccumulator chunks, occupying 32768 bytes, are not
initialized, but instead are set to zero only when actually used.  For
many applications, it will be typical for only a small fraction of the
chunks to be used, because the numbers summed have limited range, or
are all the same sign.

It is also necessary to transfer all large superaccumulator chunks to
the small superaccumulator when the final rounded result is required.
The obvious way of doing this would be to look at all 4096 counts,
transferring the corresponding chunk if the count is not $-1$. The
overhead of this can be reduced by keeping an array of 64 flag words,
each a 64-bit unsigned integer, whose bits indicate which chunks have
been used.  These flag words can be used to quickly skip large regions
of unused chunks.  This scan can be further sped up, in many cases,
using a 64-bit unsigned integer whose bits indicate which of the 64
flag words are not all zero.  Maintaining these flag words slightly
increases the cost of processing a chunk when its count is negative,
but does not increase the cost of the inner summation loop.

A naive count of operations for adding one term in the C code of
Figure~\ref{fig-large-add} gives only about 8, compared to about 19 in
Figure~\ref{fig-small-add}.  And indeed, we will see below that
summing large arrays using a large superaccumulator is about twice as
fast as summing them using only a small superaccumulator.

\subsection*{Performance evaluations}\vspace*{-7pt}

The relative performance of different method for summing the elements
in an array will depend on many factors.  Some concern the problem
instance --- such as the number of terms summed, and the range of
numerical values spanned by those terms.  Others concern the computing
environment --- such as the architecture of the processor, the speed
of memory, and which compiler is used.  There will also be random
noise in measurements.  

The resulting variability has led Langlois, Parello, Goossens, and
Porada (2012) to despair of obtaining meaningful times on real
machines, and to instead advocate assessing exact summation methods
based on reproducible measurements from a simulation of how long a
program would run on a hypothetical ideal processor in which
instruction-level parallelism (ILP) allows each operation to be
performed as soon as the operands it depends on have been computed.
While this work does provide insight into the methods they assess,
it does not answer the practical question of how well the methods
perform on real computers.

Here, I will use time measurements on real computer systems to assess
performance in a way that is both directly useful and provides
some insight into the factors affecting performance of summation
methods.  Sixteen computer systems with a variety of characteristics
were used, many of them in conjunction with several compilers.

I limited the scope of this assessment to serial implementations.
Although many of the processors used have multiple cores or threads,
only a single thread was executing during these tests.  (The systems
were largely idle apart from the test program itself.)

The small and large superaccumulator methods were implemented in C,
with careful attention to efficiency.  Several code segments were
implemented twice, once in a straightforward manner (without obvious
inefficiencies), and a second time with attempts at manual
optimizations, such as loop unrolling and branch avoidance. The
straightforward implementation might be the most efficient, if the
compiler produces superior optimization decisions.  This was not found
to be the case, however, so the manually-optimized versions were used.

The simple summation routines were similarly implemented (with
manually optimized versions chosen).  The ordered summation routine
adds each term in turn to a 64-bit double-precision accumulator.  The
unordered summation routine uses separate accumulators for terms with
even and odd indexes, then adds them together at the end.  This allows
scope for instruction-level parallelism.

For the iFastSum and OnlineExact methods of Zhu and Hayes (2010), I
used the C++ implementation provided by them as supplementary
information to their paper.  From casual perusal, this C++ code
appears to be a reasonably efficient implementation of these methods,
but it is possible that it could be improved.

The C/C++ compilers used were gcc-4.6, gcc-4.7, gcc-4.8, gcc-4.9,
clang-3.4, clang-3.5, and clang-3.6.  For many of the systems, more
than one of these compilers were available.  Choice of compiler
sometimes had a substantial impact on the performance of the various
methods, and the most recent compiler version was not always the best.
Since relative as well as absolute performance differed between
compilers, an arbitrary choice would not have been appropriate.
Instead, for each method a best choice of compiler from among those
available was made, based on the time summing 1000 terms for the small
superaccumulator and iFastSum methods, on the time summing 10000 terms
for the two simple summation methods, and on the time summing 100000
terms for the large superaccumulator and OnlineExact methods.  The
compiler chosen for each method was then used for summations
of all sizes done with that method.

Seven array sizes were tried, ranging from $N=10$ to $N=10^7$ by
powers of ten, which covers the sizes relevant to the sizes of data
caches in the processors tested.  This range of sizes also shows the
effects of fixed versus per term costs for the various methods.  Each
summation was repeated $R\, =\, 10^9 / N$ times, and the total time for
all summations was recorded, along with the total time divided by the
total number of terms summed (which was always $10^9$), which was
reported in nanoseconds per term, and is what is shown in the plots
below.

Note that due to the $R$-fold repetition, with $R$ at least 10,
summing arrays of a size for which all the data fits in the memory
cache should result in most memory accesses being to cache.  The
processors used all have at least two levels of cache, whose sizes are
shown by vertical lines in the plots, at the number of terms for which
the data would just fit in that level cache.

In order to limit the effort needed for this assessment, I mostly used
only a single distribution for numeric elements of the arrays summed,
as follows.  The terms in the first half of each array that was summed
were independent, with values given by $U_1 \exp(30U_2)$, with $U_1$
and $U_2$ being pseudo-random values uniformly drawn from $(0,1)$
using a multiplicative congruential generator with period 67101322.
(The standard C \verb|rand| generator was avoided, since it is not the
same on all systems.)  The terms in the second half of the array were
the negations of the mirror reflection of the terms in the first half
--- that is, element $N\!-\!1\!-\!i$ was the negation of element $i$.
The exact sum of all terms was therefore zero.  I also performed a few
tests in which the elements of the array were randomly permuted before
being summed, as discussed after the main results shown in the
figures.

\begin{figure}[p]

~~~~~~~~ \includegraphics[scale=0.6]{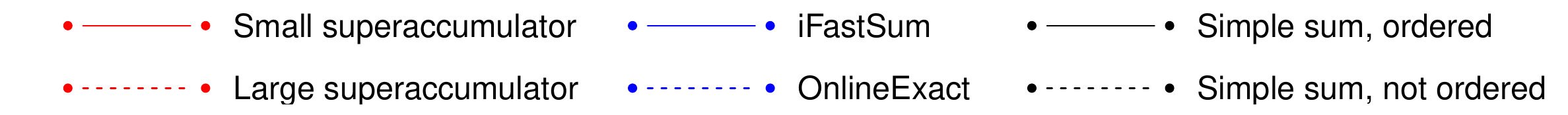}

\vspace{8pt}

~~~~~~~ \makebox[3.1in]{\small Intel Core 2 Duo (T7700), 2.4 GHz, 2007} ~~~
        \makebox[3.1in]{\small Intel Xeon E5462, 2.8 GHz, 2007}

~~ \includegraphics[scale=0.6]{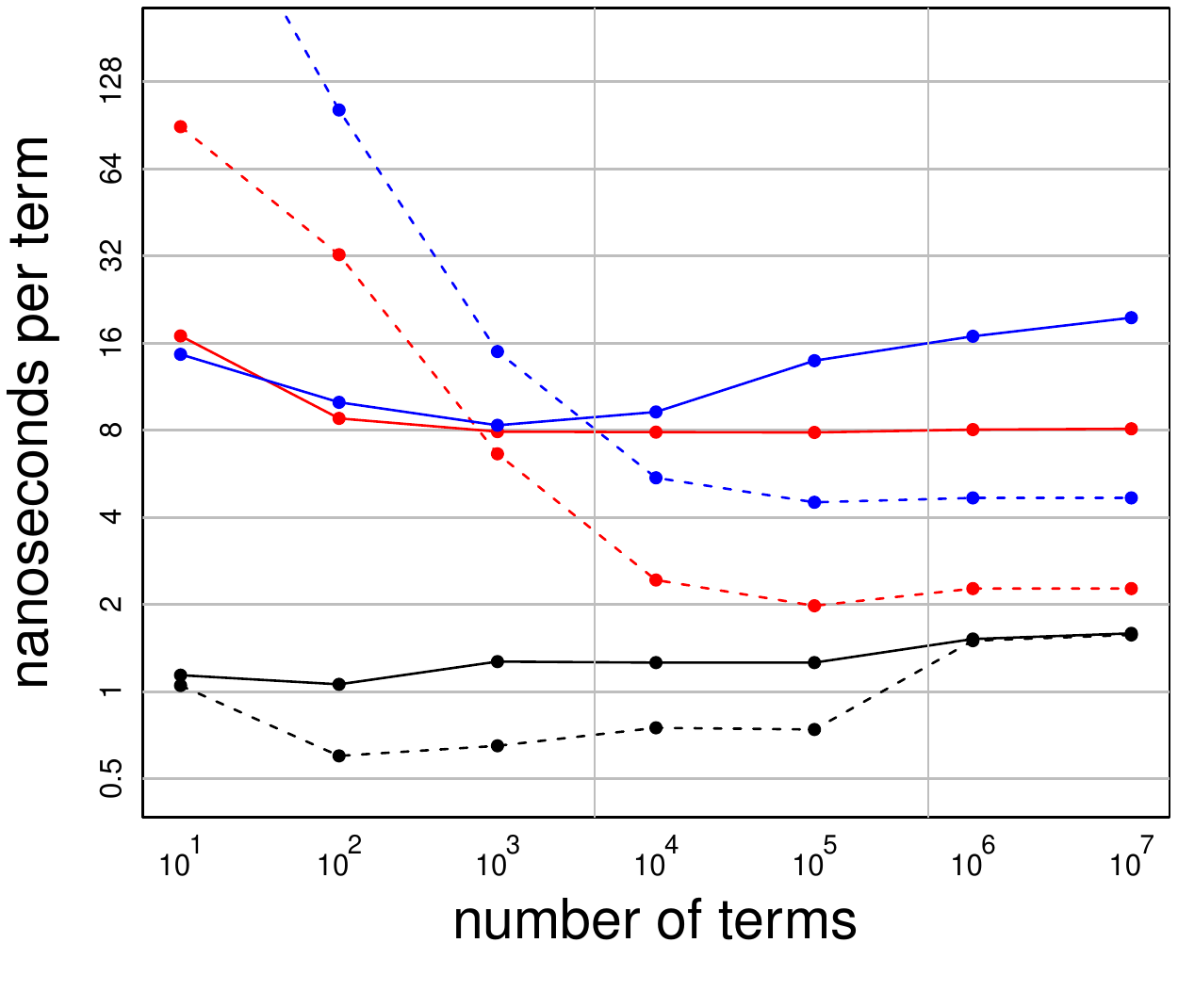} ~~~
   \includegraphics[scale=0.6]{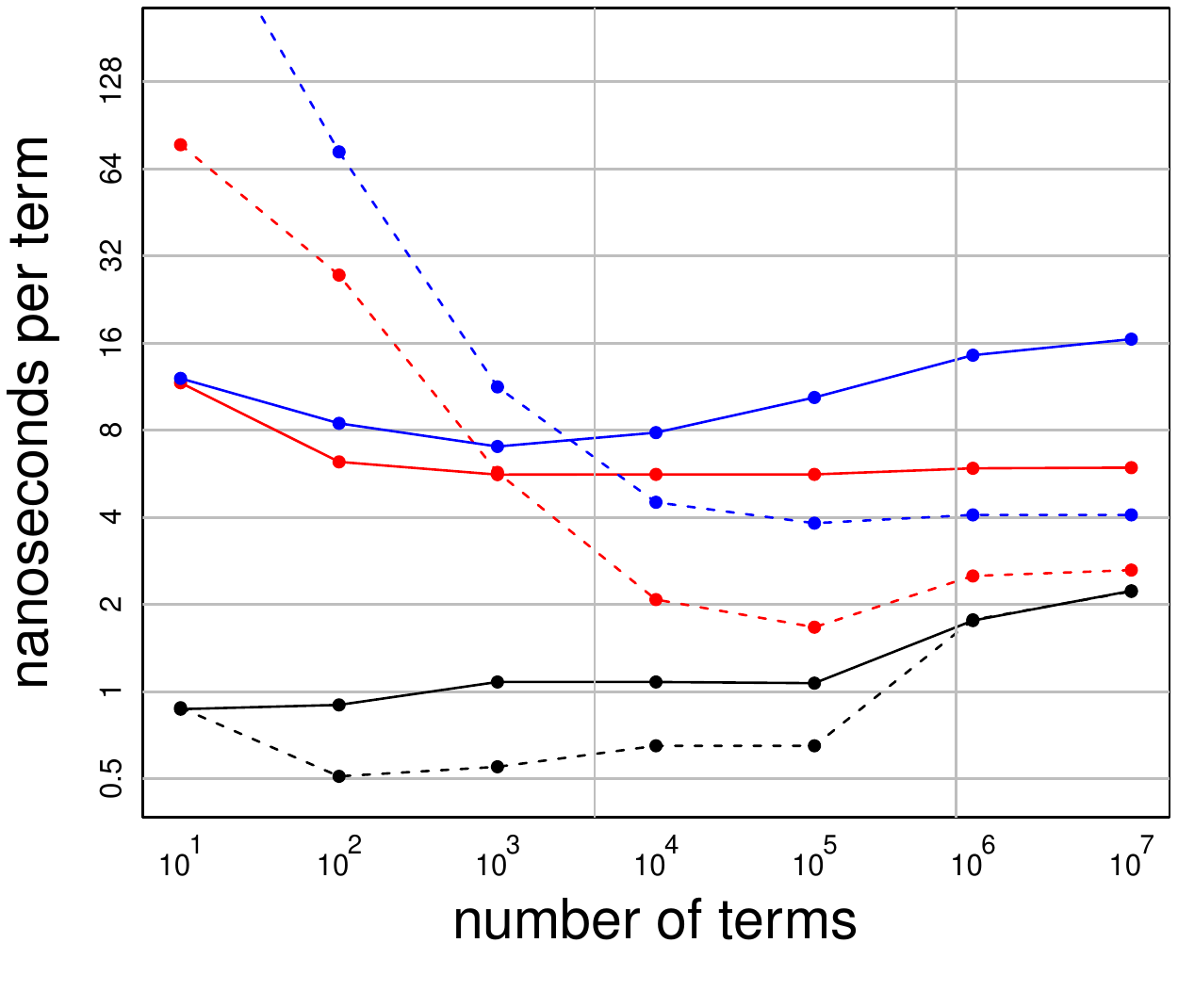}

~~~~~~~ \makebox[3.1in]{\small Intel Xeon X5680, 3.33 GHz, 2010} ~~~
        \makebox[3.1in]{\small Intel Xeon E3-1225, 3.1 GHz, 2011}

~~ \includegraphics[scale=0.6]{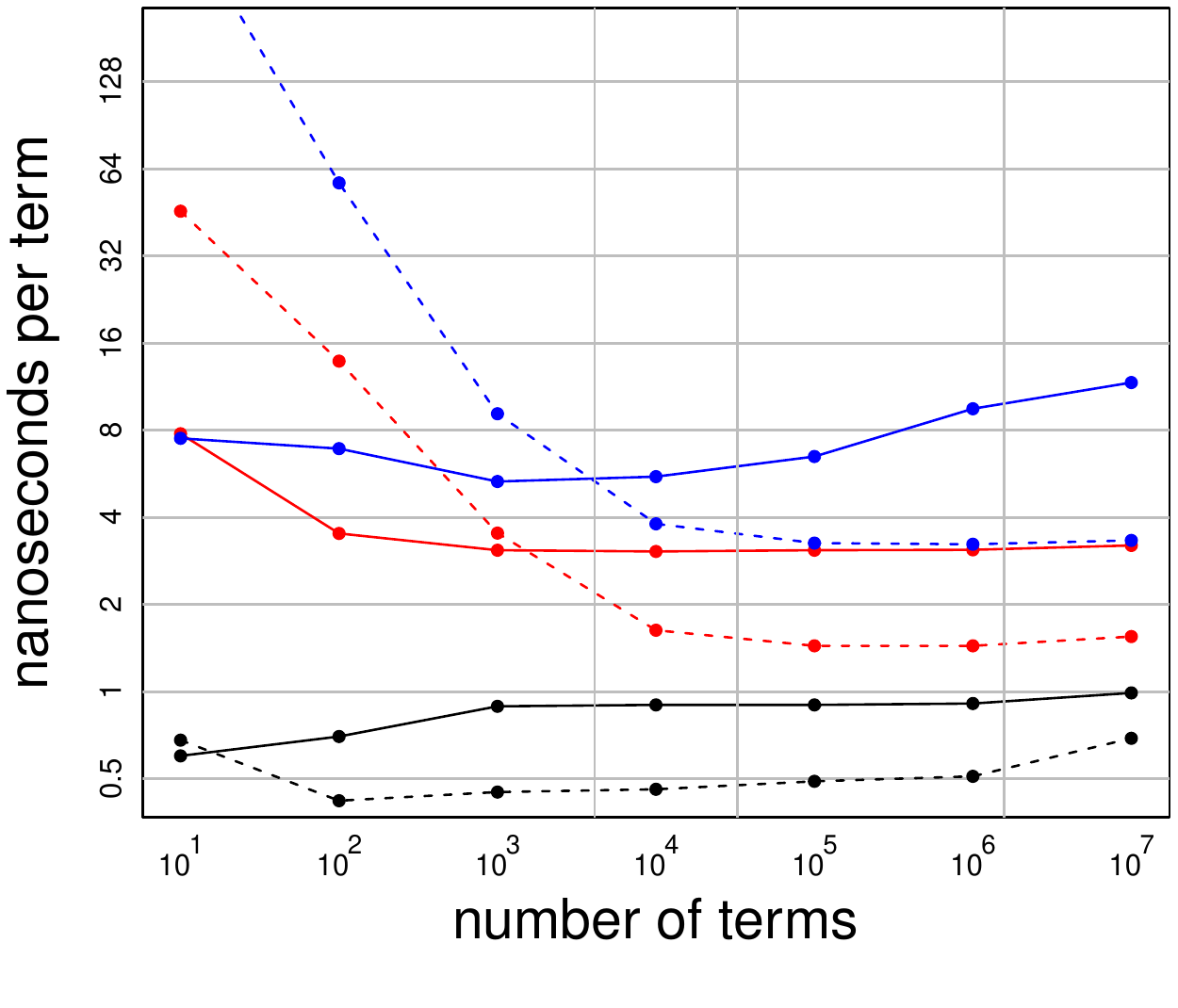} ~~~
   \includegraphics[scale=0.6]{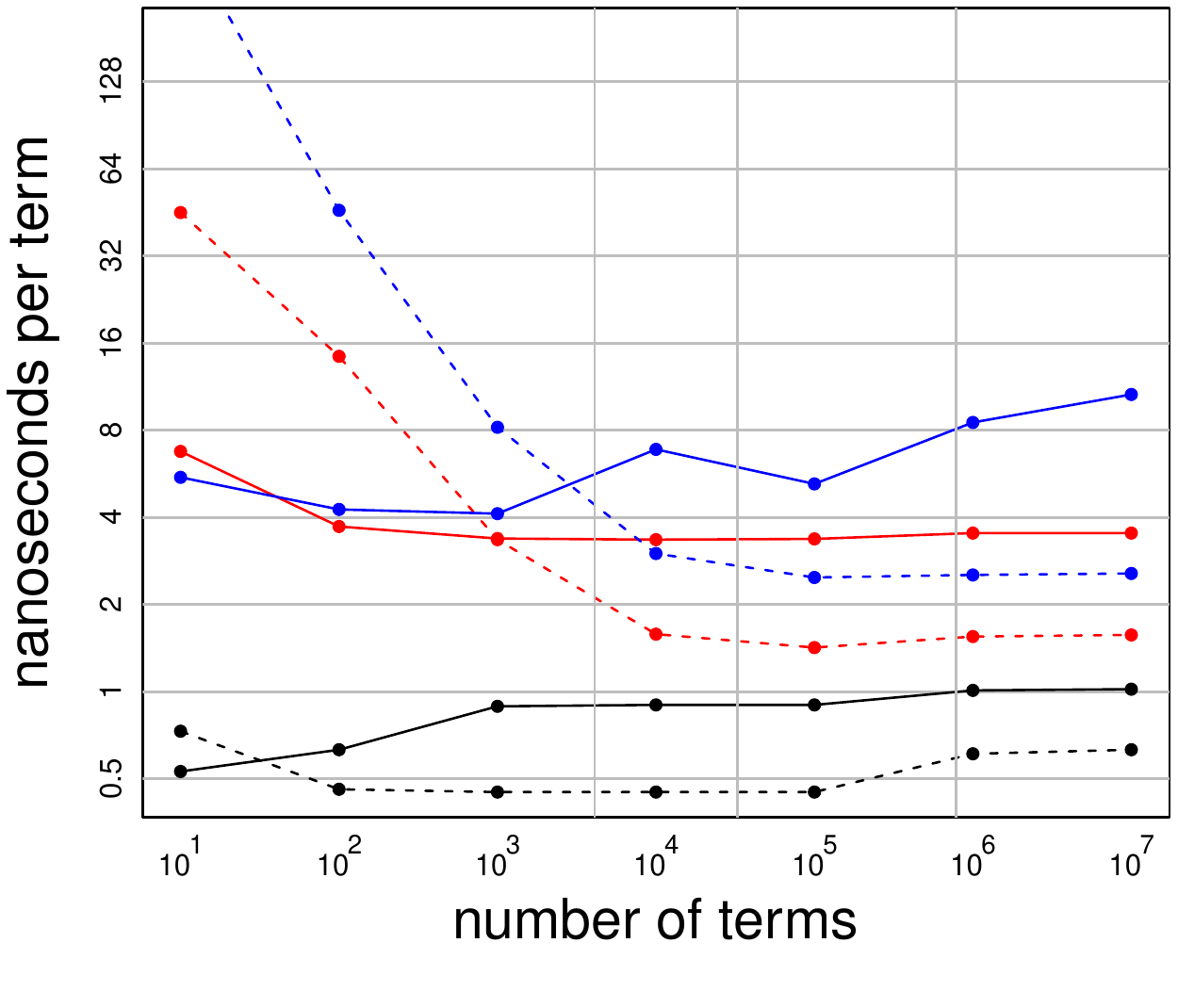}

~~~~~~~ \makebox[3.1in]{\small Intel Xeon E3-1230 v2, 3.3 GHz, 2012} ~~~
        \makebox[3.1in]{\small Intel Celeron 1019Y, 1 GHz, 2013}

~~ \includegraphics[scale=0.6]{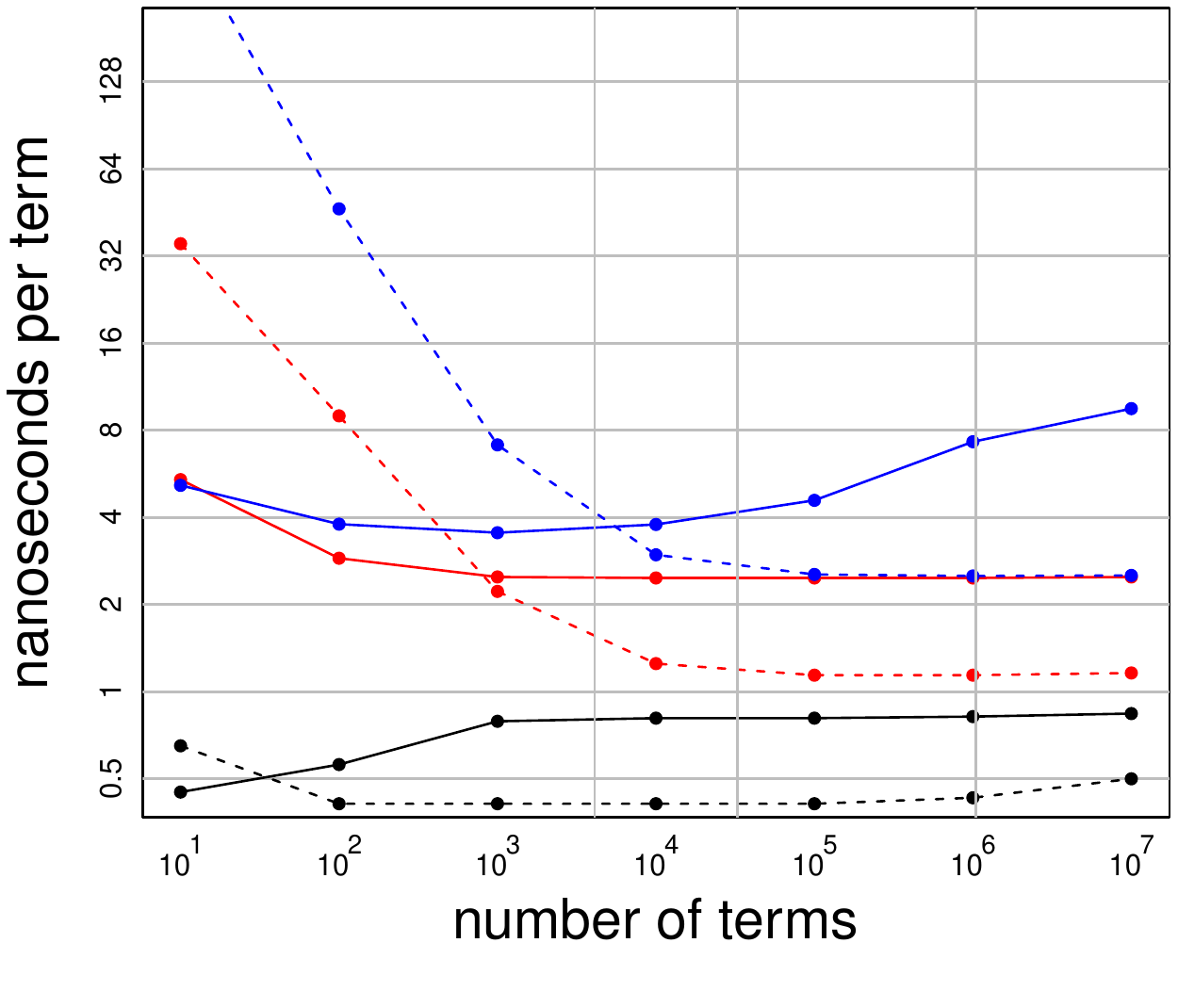} ~~~
   \includegraphics[scale=0.6]{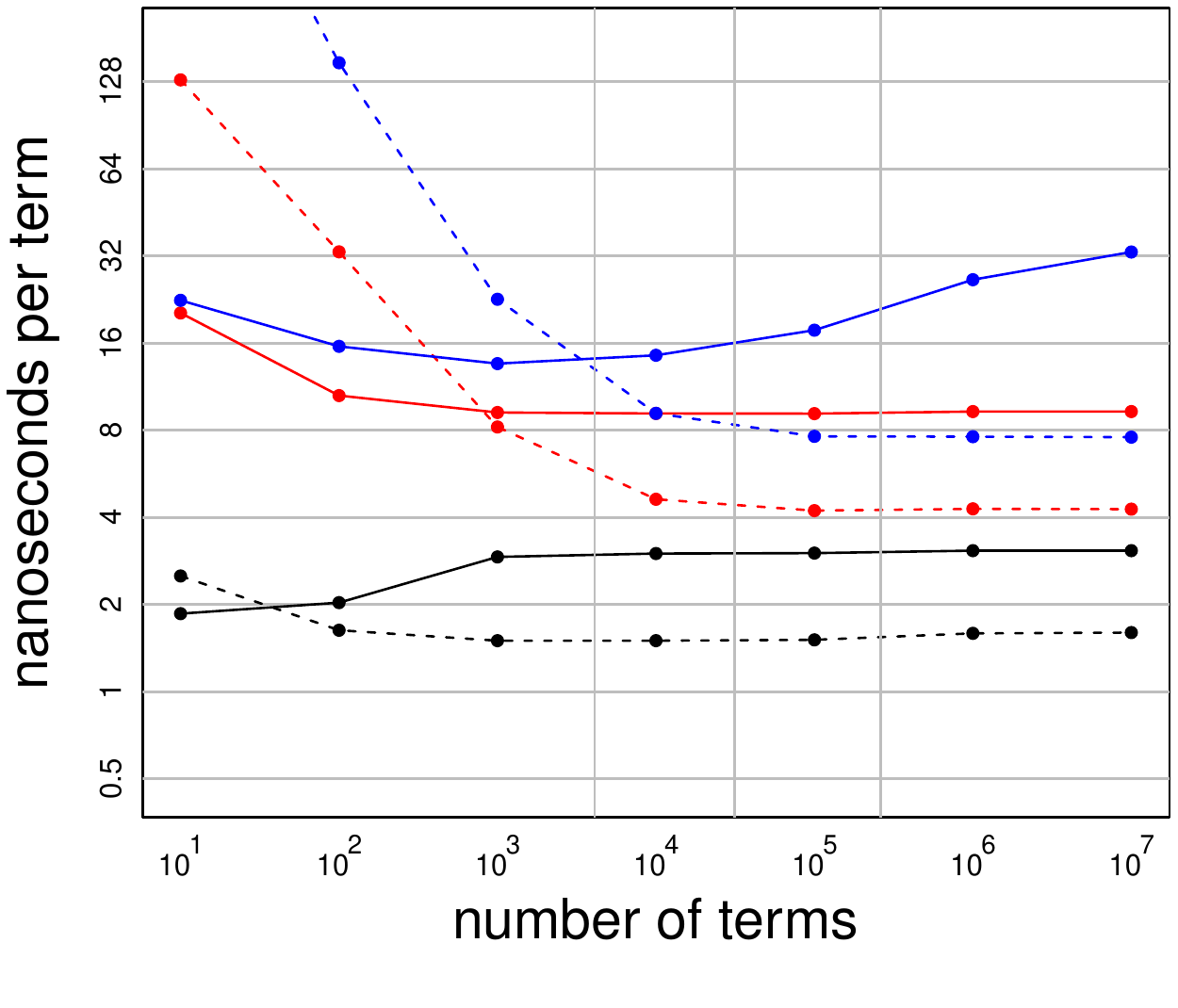}

\vspace*{-11pt}

\caption{Performance of summation methods on six 64-bit Intel 
systems.}\label{fig-res-intel64}

\end{figure}

\begin{figure}[t]

\mbox{~}

~~~~~~~ \makebox[3.1in]{\small UltraSPARC III, 1.2 GHz, 2005?} ~~~
        \makebox[3.1in]{\small UltraSPARC T2 Plus, 1.2 GHz, 2008}

~~ \includegraphics[scale=0.6]{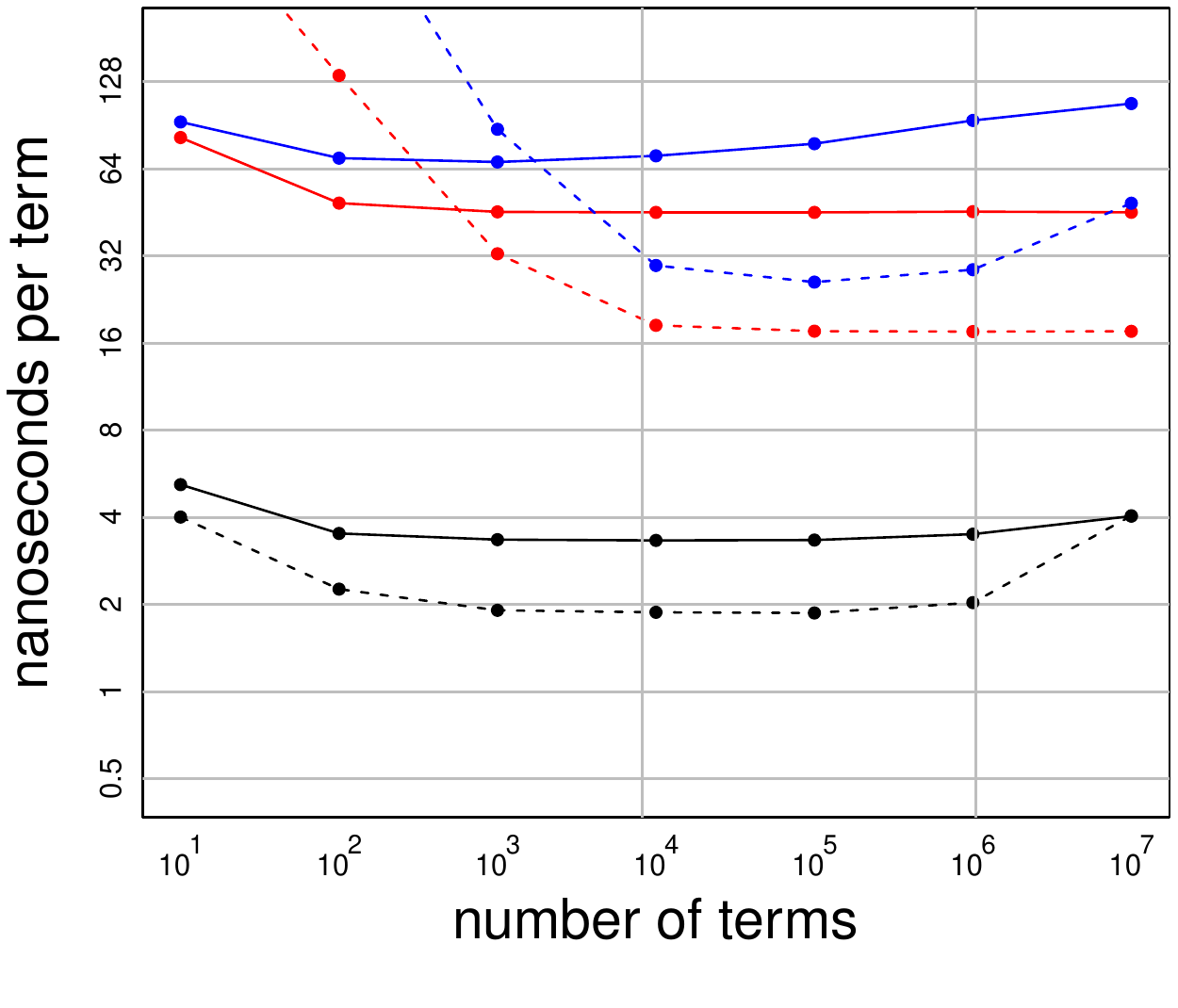} ~~~
   \includegraphics[scale=0.6]{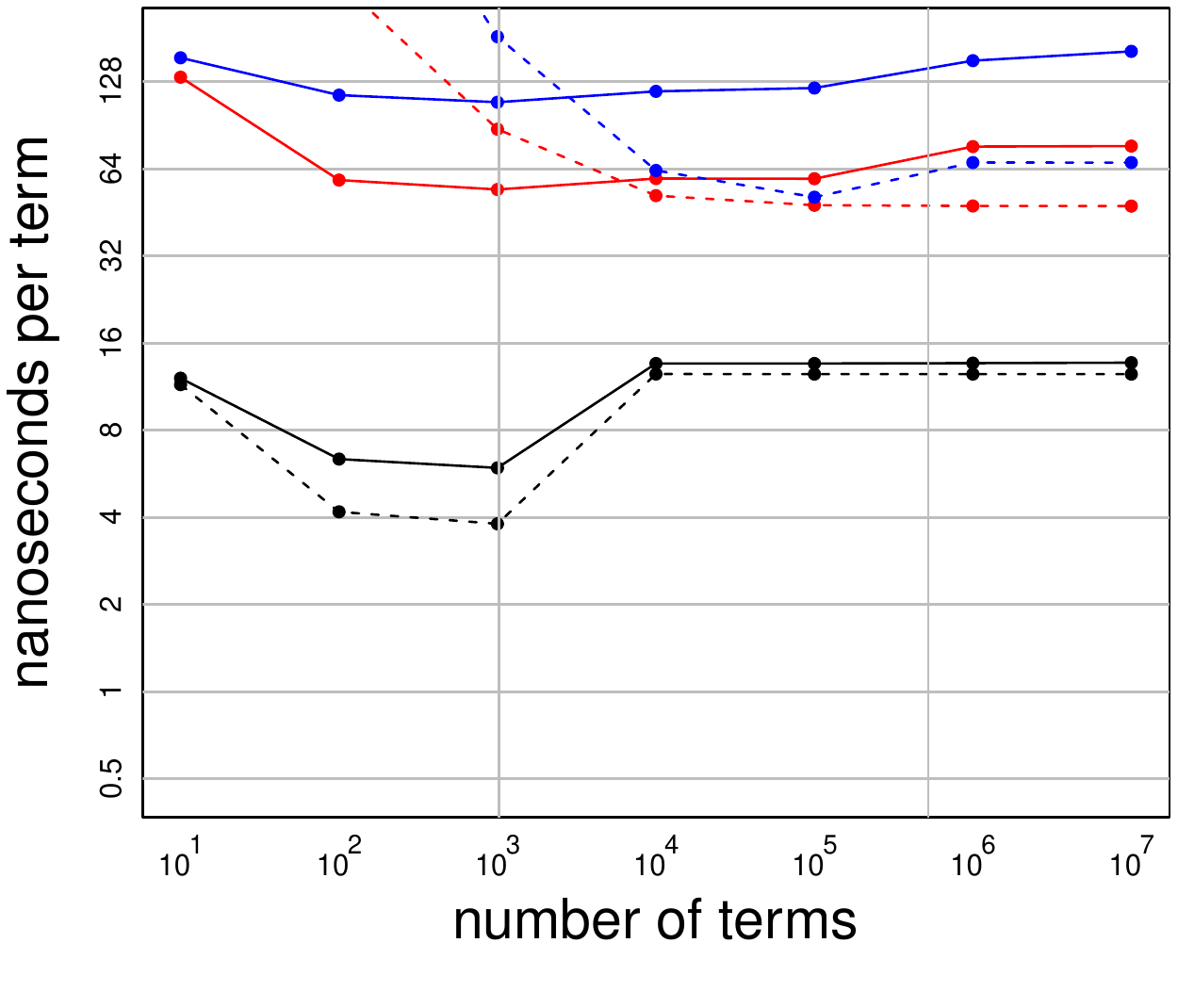}

\vspace*{-11pt}

\caption{Performance of summation methods on two 64-bit 
SPARC systems.}\label{fig-res-sparc}

\end{figure}

Figures~\ref{fig-res-intel64} though~\ref{fig-res-sparc} show the
results of the performance tests, with the six methods indicated by colour and
solid vs.\ dashed lines as shown in the key above
Figure~\ref{fig-res-intel64}.  The processor manufacturer, model, and
year of release are show above each plot.

\begin{figure}[p]

~~~~~~~ \makebox[3.1in]{\small AMD Opteron 6348, 1.4 GHz, 2012} ~~~
        \makebox[3.1in]{\small AMD E1-2500, 1.4 GHz, 2013}

~~ \includegraphics[scale=0.6]{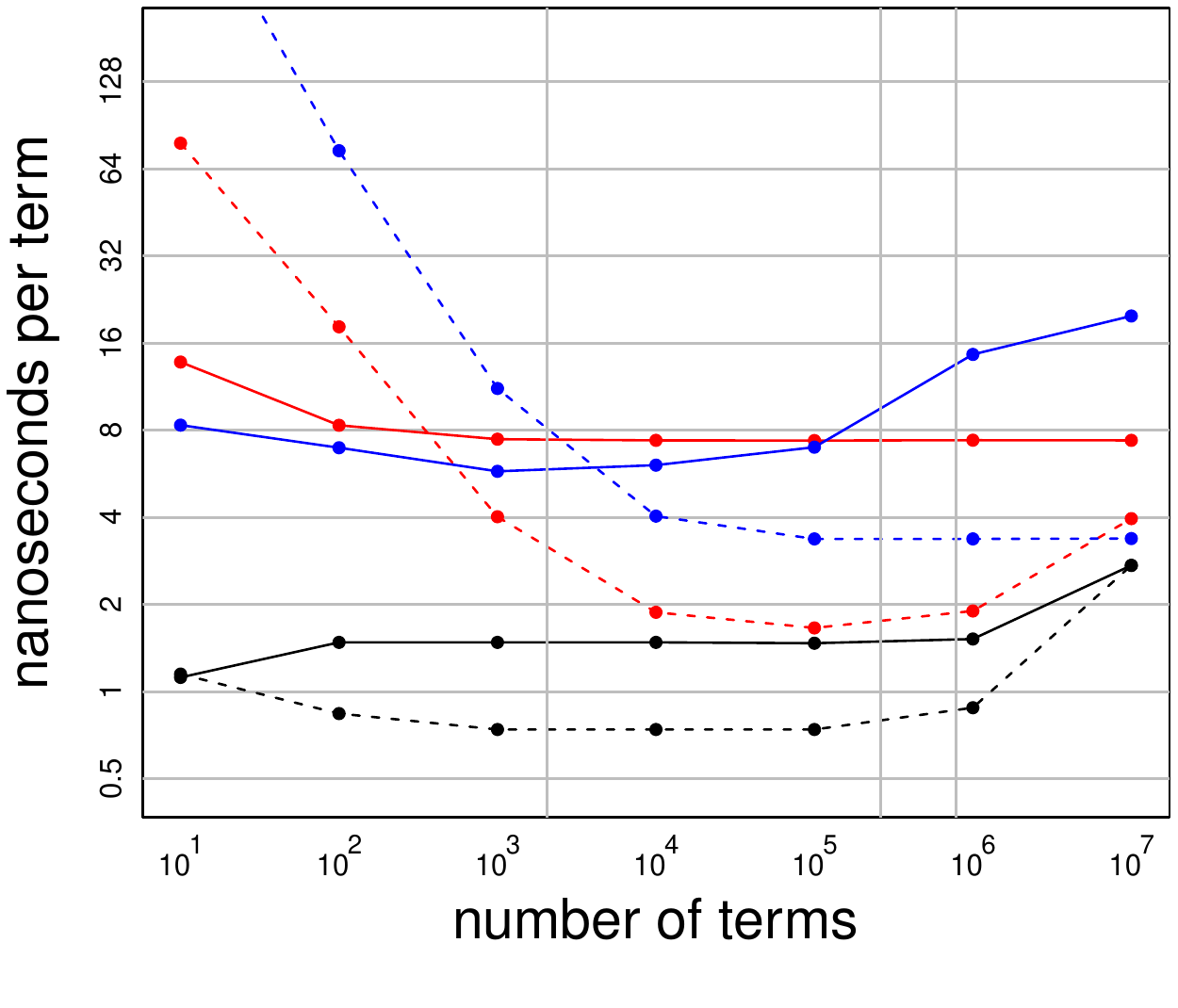} ~~~
   \includegraphics[scale=0.6]{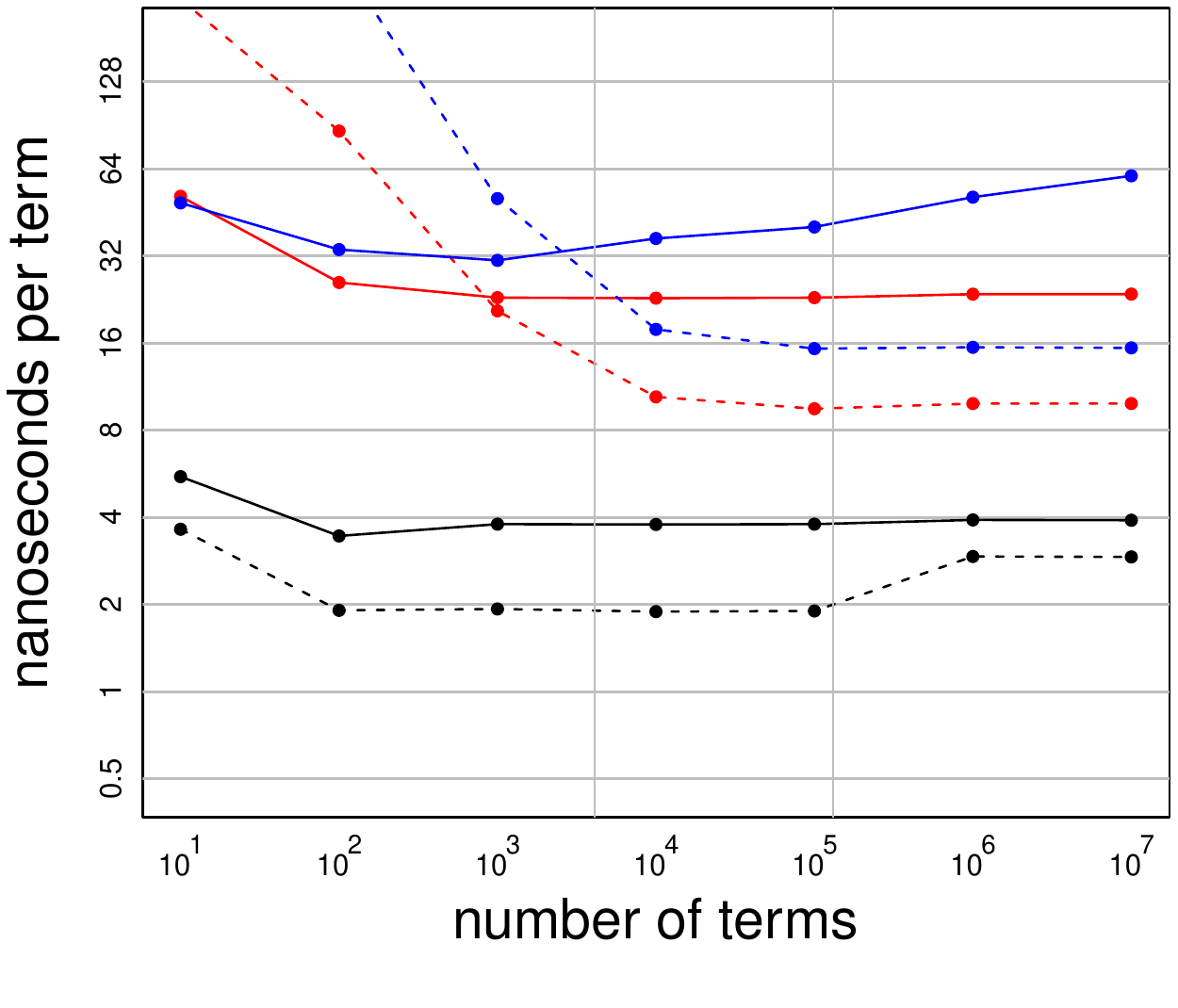}

\vspace*{-11pt}

\caption{Performance of summation methods on two 64-bit AMD 
         systems (Intel ISA).}\label{fig-res-amd}

\end{figure}

Performance on six 64-bit Intel systems and two 64-bit
AMD systems (which use the Intel Instruction Set Architecture) is
shown in Figures~\ref{fig-res-intel64} and~\ref{fig-res-amd}.  The
Xeon and Opteron processors are designed for use in servers and
high-end workstations.  The Intel Core 2 Duo is from an Apple MacBook
Pro, the Intel Celeron 1019Y is from a low-end Acer AspireV5 laptop,
and the AMD E1-2500 is from a low-end Gateway desktop system.  The six
Intel processors span three major microarchitecture families ---
``Core'' (Core 2 Duo, Xeon E5462), ``Nehalem'' (X5680), and
``Sandy/Ivy Bridge'' (Xeon E3-1225, Xeon E3-1230 v2, and Celeron
1019Y).  The two AMD processors also have different microarchitectures
--- ``Piledriver'' (Opteron 6348) and ``Jaguar'' \mbox{(E1-2500)}.


\begin{figure}[p]

\mbox{~}

~~~~~~~ \makebox[3.1in]{\small Intel Pentium III, 1 GHz, 2000} ~~~
        \makebox[3.1in]{\small Intel Xeon, 1.7 GHz, 2001}

~~ \includegraphics[scale=0.6]{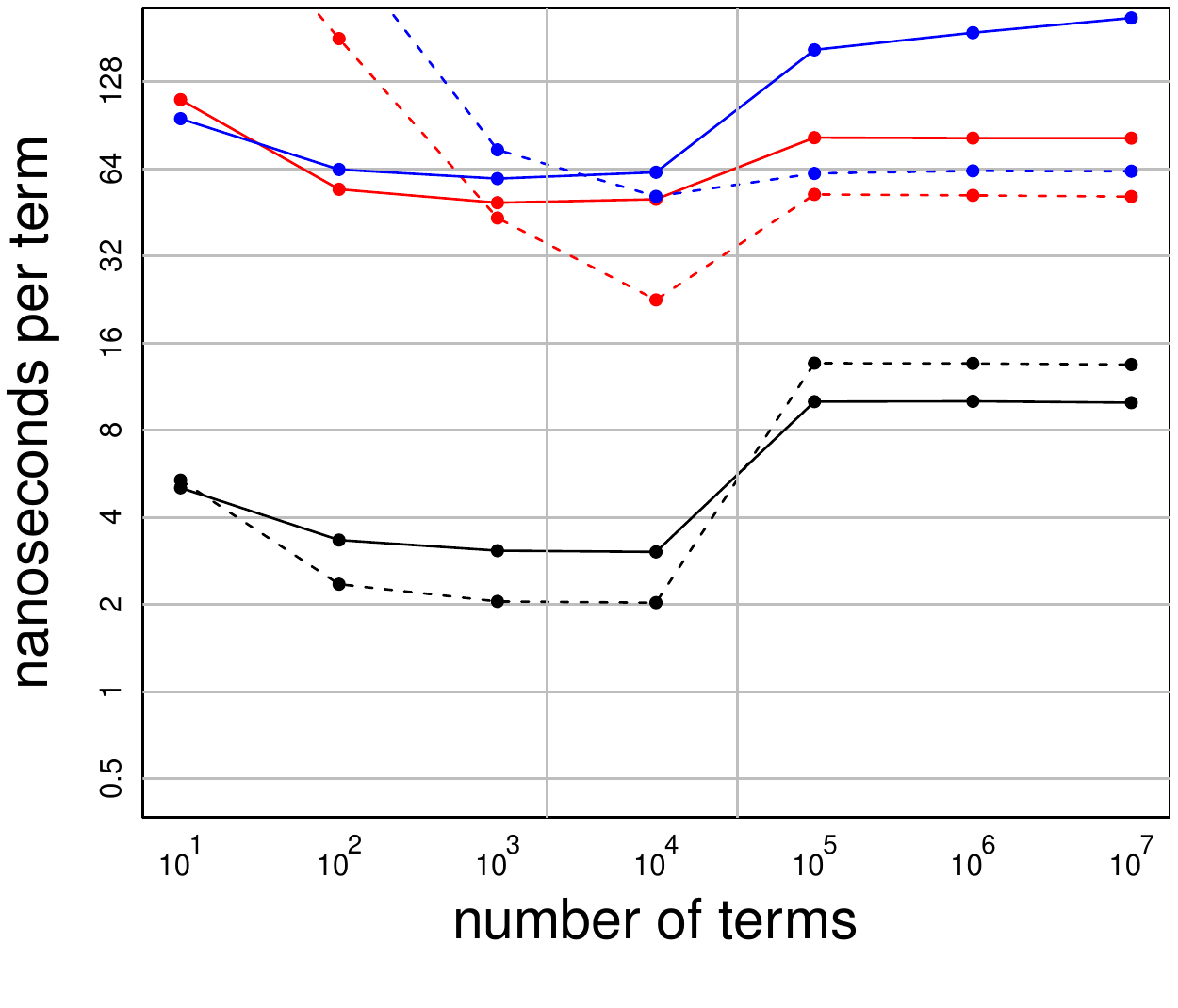} ~~~
   \includegraphics[scale=0.6]{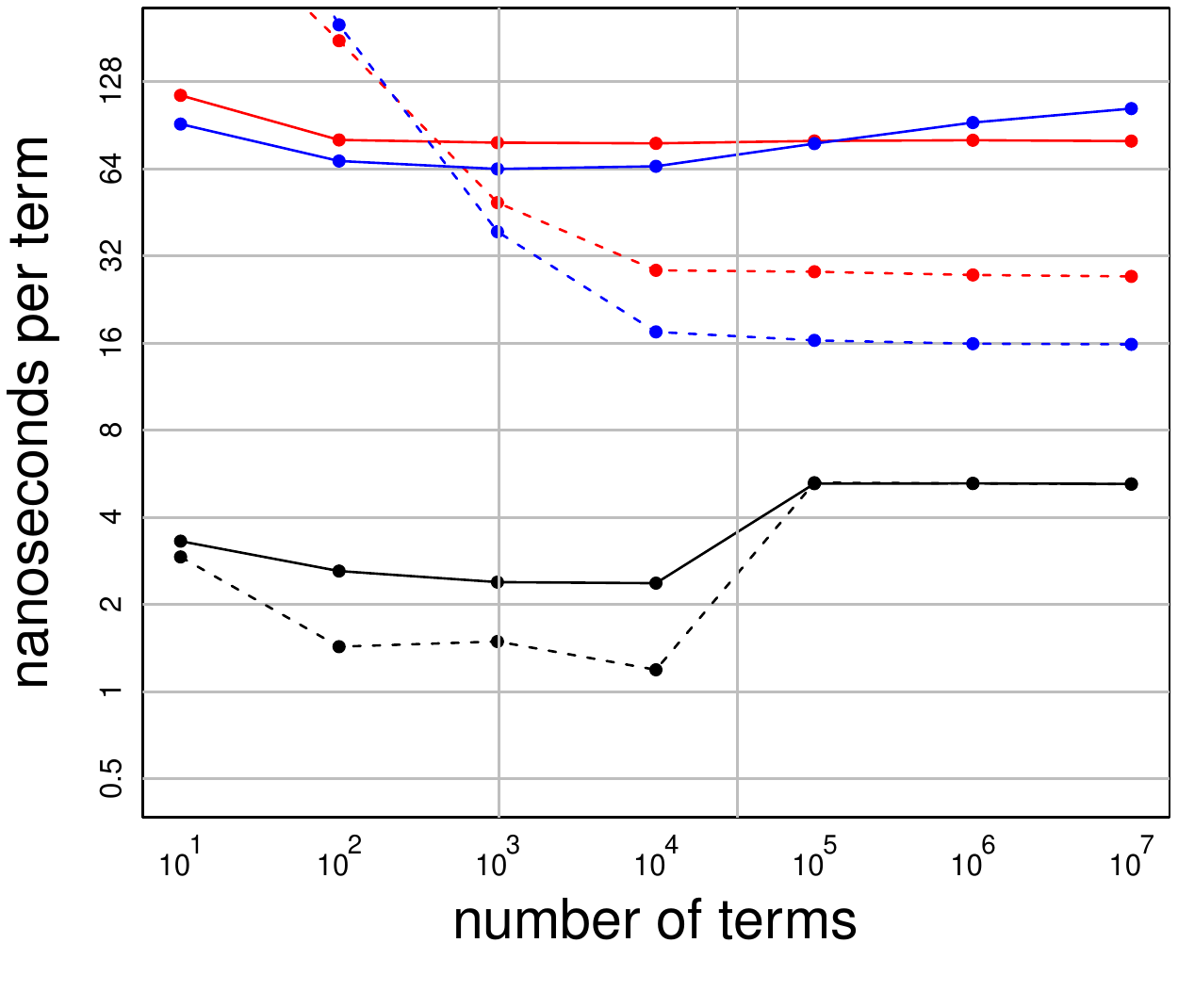}

~~~~~~~ \makebox[3.1in]{\small Intel Pentium 4, 3.06 GHz, 2002} ~~~
       \makebox[3.1in]{\small Intel Xeon X5355, 2.66 GHz, 2006}

~~ \includegraphics[scale=0.6]{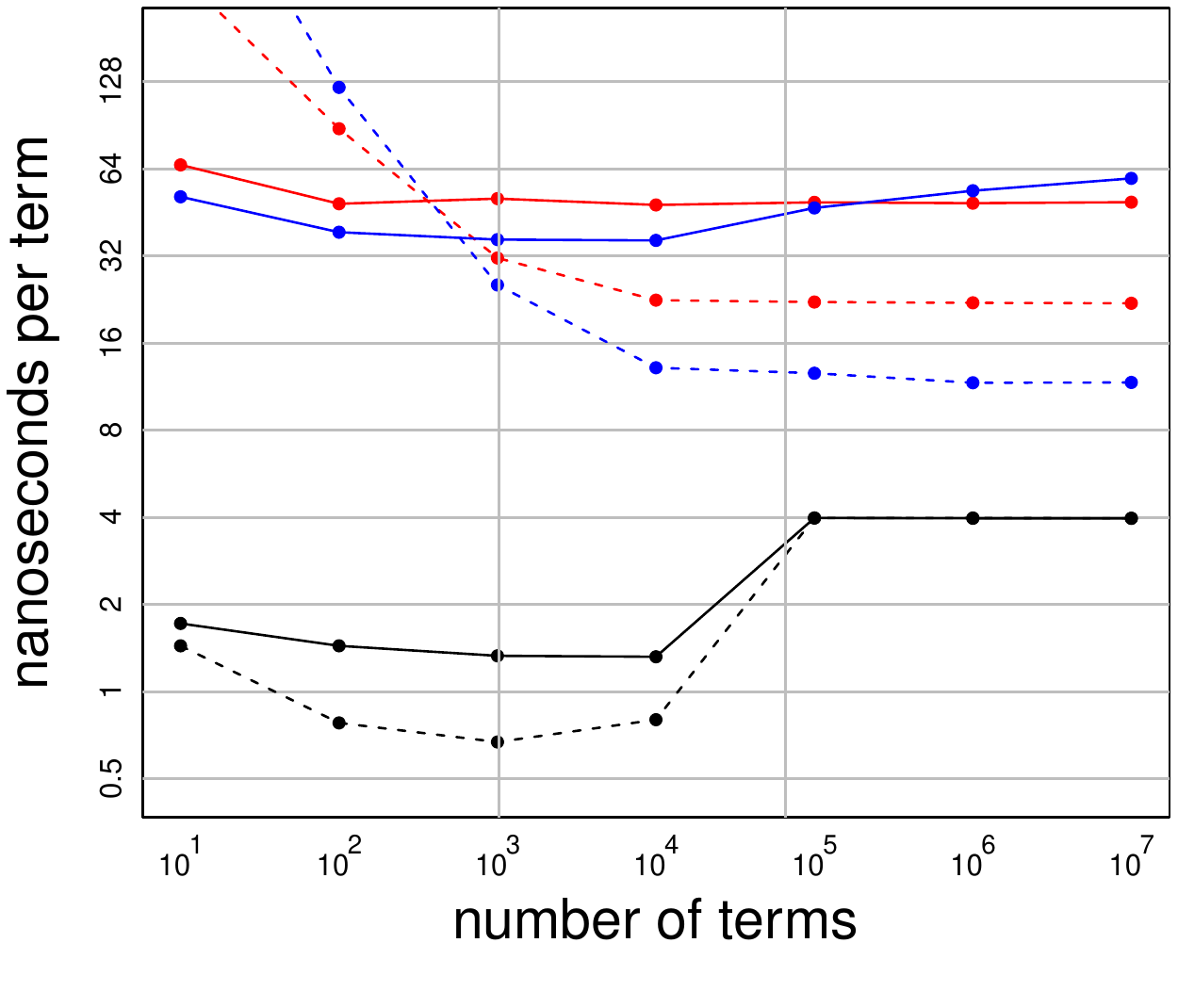} ~~~
   \includegraphics[scale=0.6]{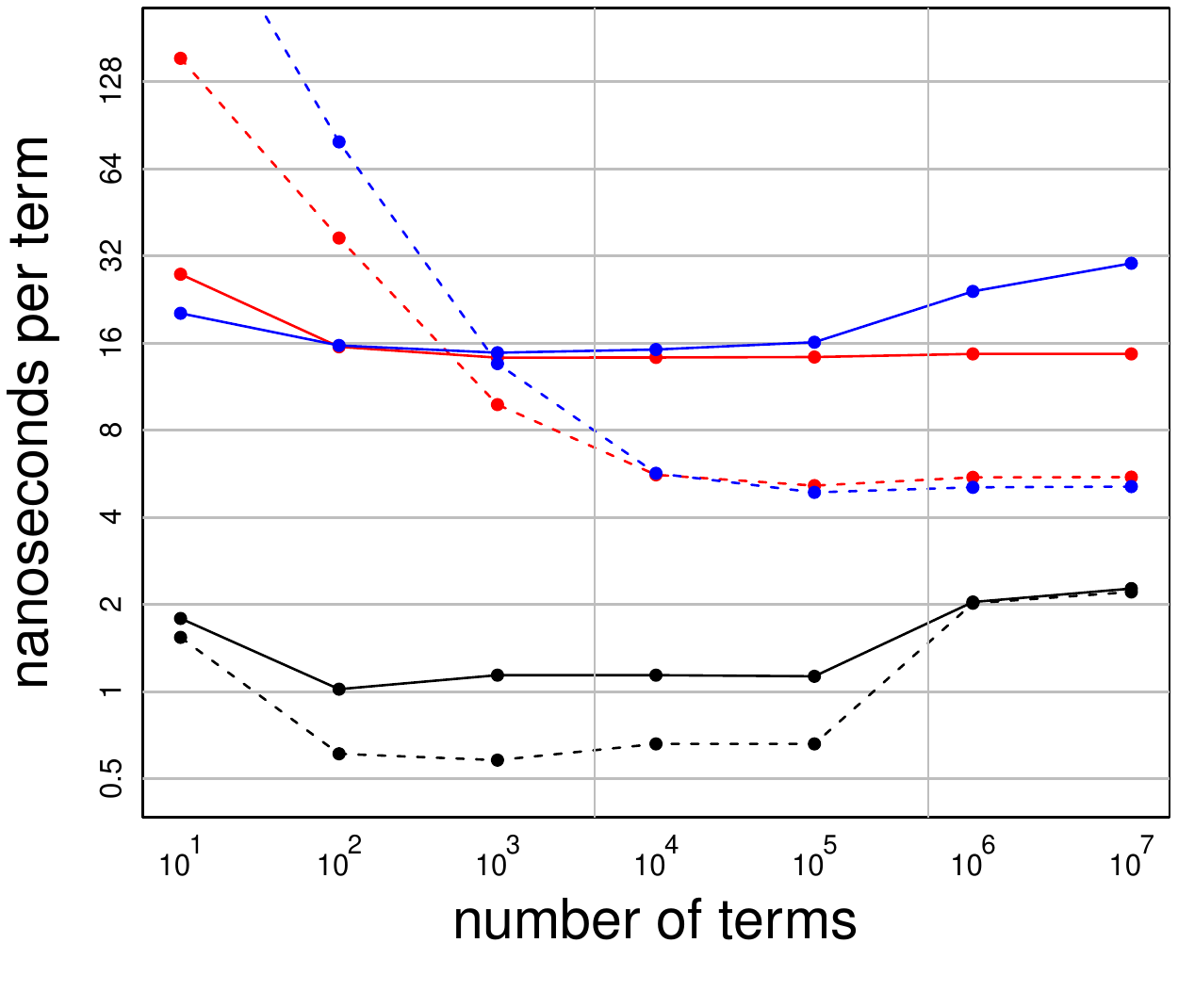}

\vspace*{-11pt}

\caption{Performance of summation methods on four 32-bit Intel systems. 
The Intel Xeon X5355 is a 64-bit capable processor, but was run in 32-bit 
mode.}\label{fig-res-intel32}

\end{figure}


\begin{figure}[t]

~~~~~~~ \makebox[3.1in]{\small ARMv6 Processor, 700 MHz, 2003} ~~~
        \makebox[3.1in]{\small Cortex-A9 ARMv7 Processor, 1 GHz, 2008}

~~ \includegraphics[scale=0.6]{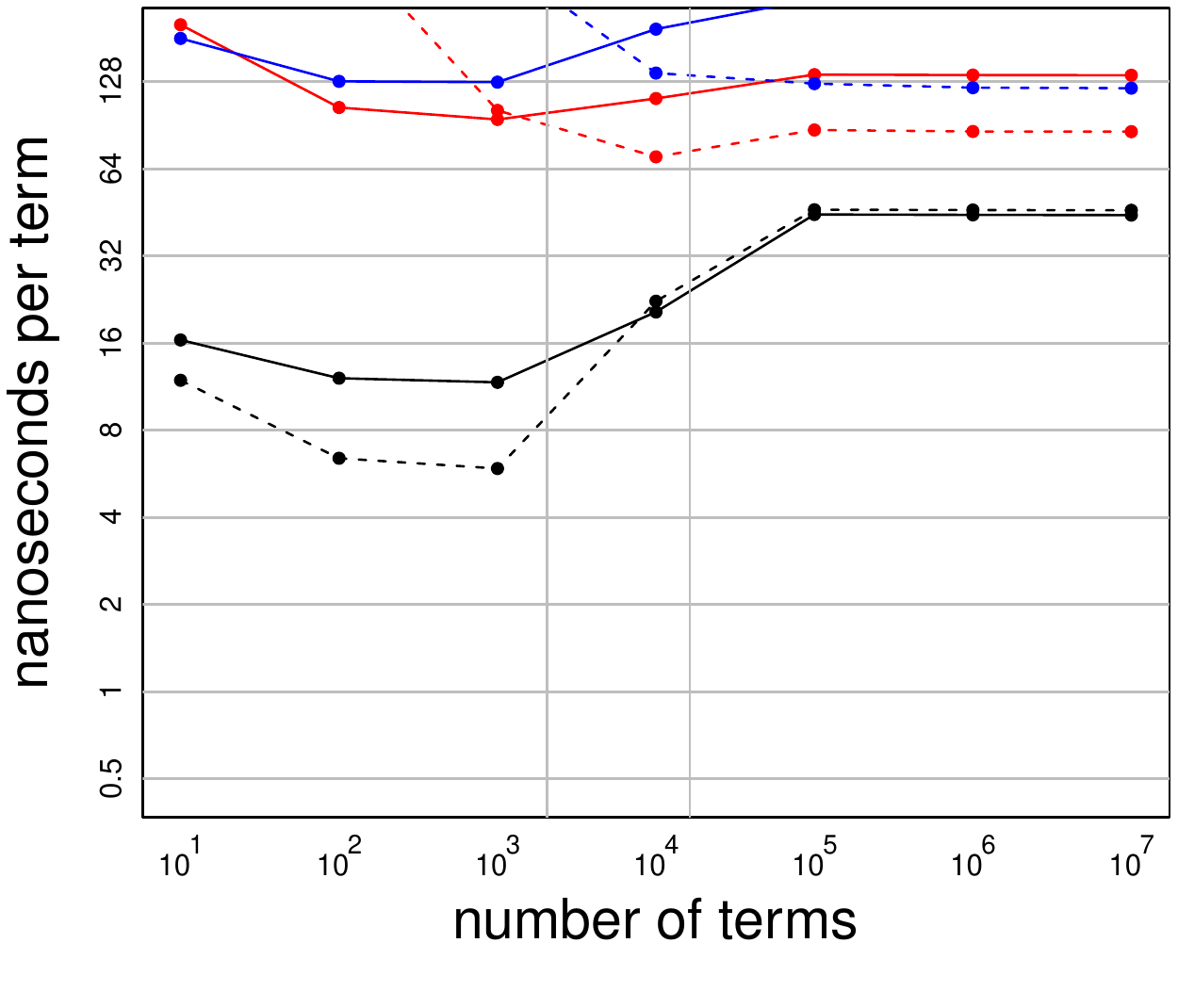} ~~~
   \includegraphics[scale=0.6]{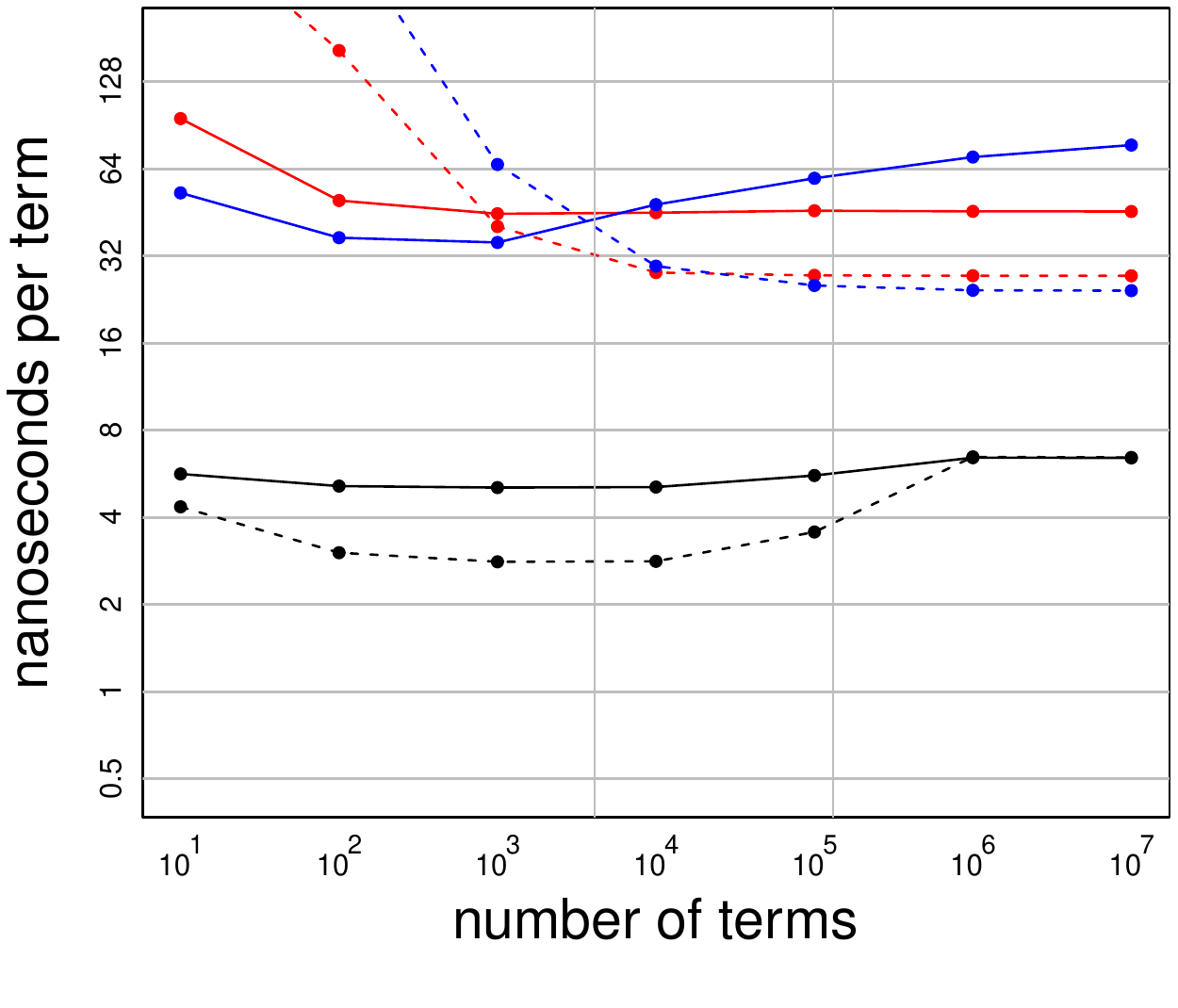}

\vspace*{-11pt}

\caption{Performance of summation methods on two 32-bit ARM 
systems.}\label{fig-res-arm}

\end{figure}


The qualitative picture from these tests on modern processors is quite
consistent.  The large superaccumulator method is faster than the
small superaccumulator method when summing more than about 1000 terms.
Similarly, the OnlineExact method is faster than the iFastSum method
when summing more than about 2500 terms.  The combination of the two
superaccumulator methods --- the small superaccumulator method for
less than 1000 terms, and the large superaccumulator method for 1000
terms or more --- is superior to any combination of the iFastSum and
OnlineExact methods, except that for some processors iFastSum is
slightly faster when summing very small arrays (less than about thirty
terms, or less a few hundred terms for the AMD Opteron 6348 processor).

The advantage of the large superaccumulator method over the
OnlineExact method for summing a large number of terms (10000 or more)
is about a factor of two, except that for some processors this
decreases (to nothing for the AMD Opteron 6348) when summing very
large arrays, for which out-of-cache memory access time dominates.
The advantage of the small superaccumulator method over iFastSum when
summing small arrays is less (non-existent for the Intel Core 2 Duo
and the AMD Opteron 6348), but the small superaccumulator method
nevertheless appears to be generally preferable to iFastSum for other
than very small sums.

The large superaccumulator method is no more than about a factor of
two slower than simple ordered summation, when summing 10000 or more
terms.  For the AMD Opteron 6348, the large superaccumulator method is
only slightly slower than simple summation, though for the AMD
\mbox{E1-2500} the ratio of times is slightly greater than two.  The
difference between the large superaccumulator method and simple ordered
summation is often less for very large summations, as expected if time
for out-of-cache memory accesses starts to dominate.

The simple summation method that adds terms out of order is about
twice as fast as simple ordered summation, except for summing very
large arrays, for which its advantage is usually less (sometimes
non-existent).  For array sizes of 10000 and 100000, simple unordered
summation is typically about three times faster than the large
superaccumulator method.

How the methods perform on four 32-bit Intel processors is shown in
Figure~\ref{fig-res-intel32}.  The Intel Pentium III processor uses
the ``P6'' microarchitecture, which is a distant ancestor of the
``Core'' microarchitecture of the Intel Xeon X5355.  The Intel Xeon
(1.7 GHz) and Intel Pentium 4 use the ``NetBurst'' microarchitecture.
The Intel Pentium III processor uses the 387 floating-point unit for
floating-point arithmetic, whereas the other processors have the SSE2
floating-point instructions, which have more potential for
instruction-level parallelism.

As was the case for the 64-bit processors, we see that for summing
large arrays, the large superaccumulator method is better than the
small superaccumulator method, and the OnlineExact method is better
than the iFastSum method.  The combination of small and large
superaccumulator methods is better than the combination of iFastSum
and OnlineExact for the Intel Pentium III, except for very small
arrays.  The advantage of the large superaccumulator method over
OnlineExact is a factor of two when summing 10000 terms, but is not as
large for 1000 terms (probably because of fixed overhead) or for
100000 terms (probably because out-of-cache memory access time starts
to dominate).  For the Intel Xeon X5355, the large superaccumulator
method has only a slight advantage when summing 1000 terms, and the
superaccumulator methods perform almost identically to
iFastSum+OnlineExact for other sizes (with the small superaccumulator
method being slower than iFastSum for very small sums).

For the two processors with ``NetBurst'' microarchitecture -- the
Intel Xeon at 1.7 GHz and the Intel Pentium 4 --- the picture is quite
different.  For these processors, the combination of iFastSum and
OnlineExact is better than the combination of the small and large
superaccumulator methods for all array sizes.  The advantage of
OnlineExact over the large superaccumulator method is almost a factor
of two for summing large arrays.  For smaller arrays, there is less
difference between the methods.  One might speculate that this
reflects a design emphasis on floating-point rather than integer
performance in the ``NetBurst'' processors.  One can see that for
summing 10000 terms, both the small and large superaccumulator methods
are actually slower on the 1.7 GHz Xeon than on the 1 GHz Pentium III,
whereas both simple summation and the OnlineExact method perform
substantially better.

For all these 32-bit Intel processors, the ratios of the times for the
exact summation methods to the times for simple summation are
substantially greater than for the 64-bit processors (though less so
for very large summations, where out-of-cache memory access time is
large).  This also may reflect a somewhat specialized design
philosophy for these processors, in which general-purpose computation
was supported only with 32-bit registers and operations, whereas
support for 64-bit floating-point computations was similar to that
found in modern 64-bit processors.

Figure~\ref{fig-res-amd} shows result on two 32-bit ARM processors.
For the ARMv6 processor, the combination of small and large
superaccumulator methods performs better than iFastSum+OnlineExact,
but for the Cortex-A9 ARMv7 processor the comparison is mixed.
The exact summation methods are again slower compared to simple
summation than is the case for the modern 64-bit processors.

Finally, Figure~\ref{fig-res-sparc} shows results for two UltraSPARC
64-bit processors.  For both processors, the combination of the small
and large superaccumulator methods performs significantly better than
iFastSum+OnlineExact.  The performance of the superaccumulator methods
is slower compared to simple summation for these processors than for
the 64-bit Intel and AMD processors.  One should note that the
UltraSPARC T2 Plus is optimized for multi-threaded workloads, with
8~threads per core, so a performance comparison using a single thread,
as here, may be misleading.

One can measure the fixed overhead of the small superaccumulator
method by looking at the ratio of the time per term for 10 terms and
for 100 terms.  This ratio is roughly 2 for most of the processors
tested.  Assuming that the time for a summation can be modelled as
$a+bN$, where $N$ is the number of terms, $a$ is the fixed cost, and
$b$ is the per term cost, one can work out that $a/b$ is about 12.5,
as was mentioned earlier.  This model does not work well for the large
superaccumulator method, perhaps because the ``fixed'' overhead is not
actually fixed when $N$ is small, since the number of chunks used in
the large superaccumulator will grow at a substantial rate with the
number of terms summed when $N$ is still fairly small.

I did a few timing runs in which the elements of the arrays were
randomly permuted before summing them.  This has the effect of mixing
positive and negative terms randomly (rather than all positive terms
coming before all negative terms), and also affects the contents of
the superaccumulators at intermediate stages.  This permutation had
little effect on the performance of the large superaccumulator and
OnlineExact methods.  However, it did increase the time for the small
superaccumulator method on some processors, including recent ones.
This is perhaps due to the conditional branch in the inner loop shown
in Figure~\ref{fig-small-add}, which can be well-predicted if many
terms in a row have the same sign, but not if positive and negative
terms are mixed, which will affect the time on processors that do
speculative execution of instructions that may follow a branch.  The
time for the iFastSum method was not affected as much, so the relative
advantage of the small superaccumulator method over iFastSum was
smaller, though the small superaccumulator method was still faster
when summing at least 100 terms.

Kahan's (1965) method for reducing summation error (but without
producing the exact result) was also tested on all systems. On modern
64-bit processors, computing the exact sum with the large
superaccumulator method was faster than Kahan's method for summations
of more than about 1000 terms.  Kahan's method was significantly
faster than the small superaccumulator method only for summations of
less than about 100 terms.

\begin{figure}[t]

~~~~~~~ \makebox[3.1in]{\small Intel E3-1230 v2, 3.3 GHz, 2012} ~~~
        \makebox[3.1in]{\small AMD Opteron 6348, 1.4 GHz, 2012}

~~ \includegraphics[scale=0.6]{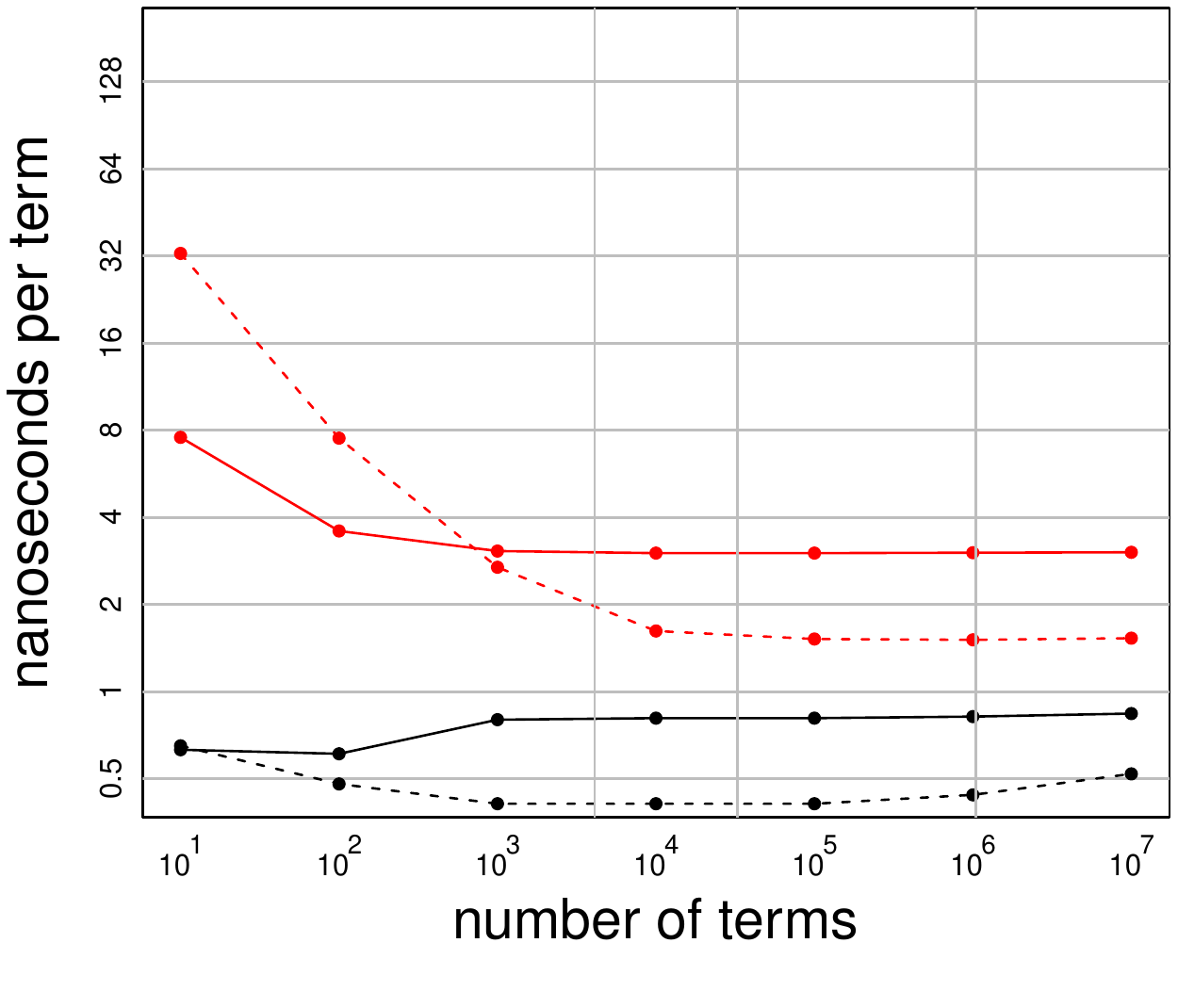} ~~~
   \includegraphics[scale=0.6]{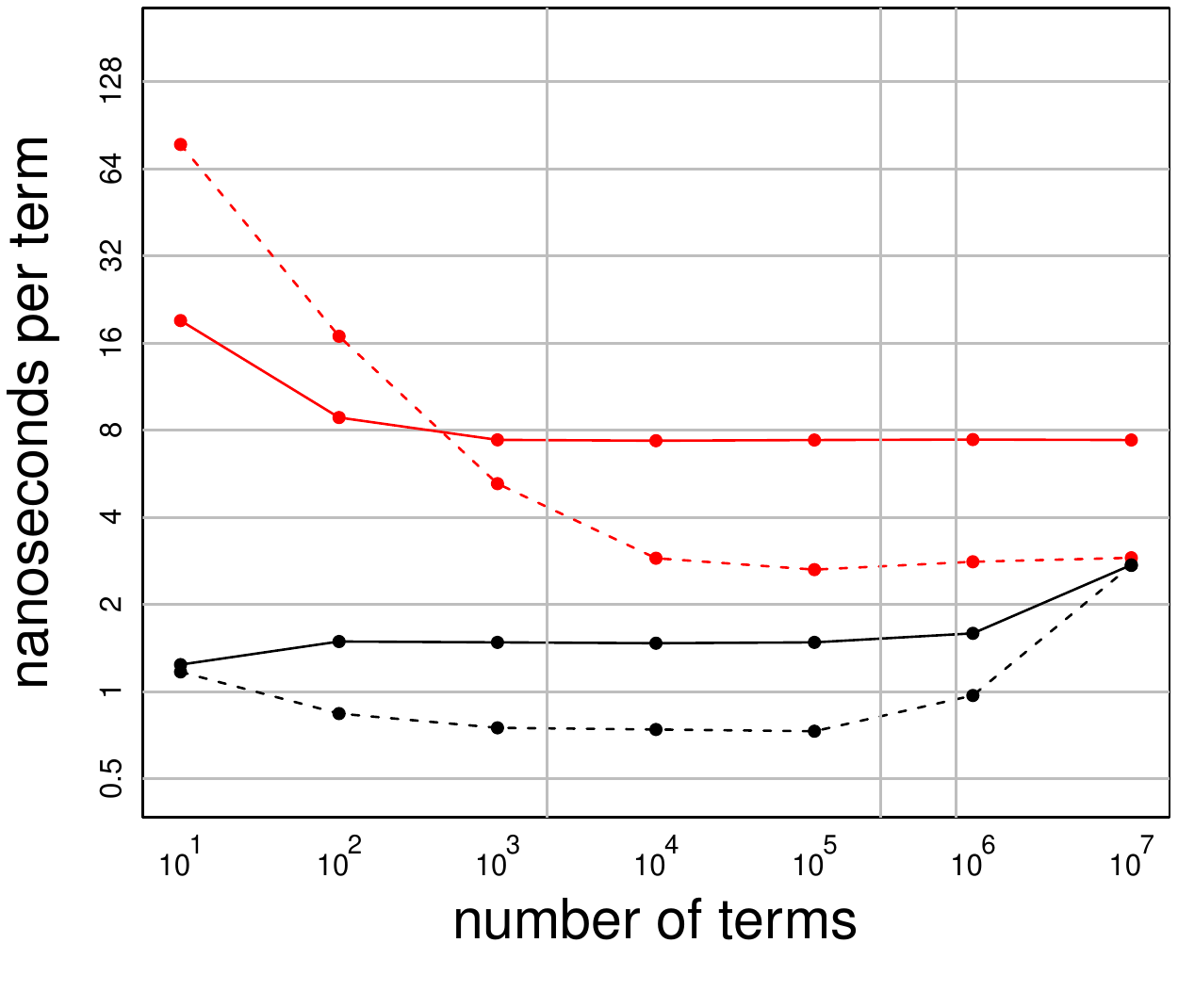}

\vspace*{-11pt}

\caption{Performance of squared norm on recent Intel and AMD high-end 
         processors.}\label{fig-res-norm}

\end{figure}

\begin{figure}[t]

\mbox{~}

~~~~~~~ \makebox[3.1in]{\small Intel E3-1230 v2, 3.3 GHz, 2012} ~~~
        \makebox[3.1in]{\small AMD Opteron 6348, 1.4 GHz, 2012}

~~ \includegraphics[scale=0.6]{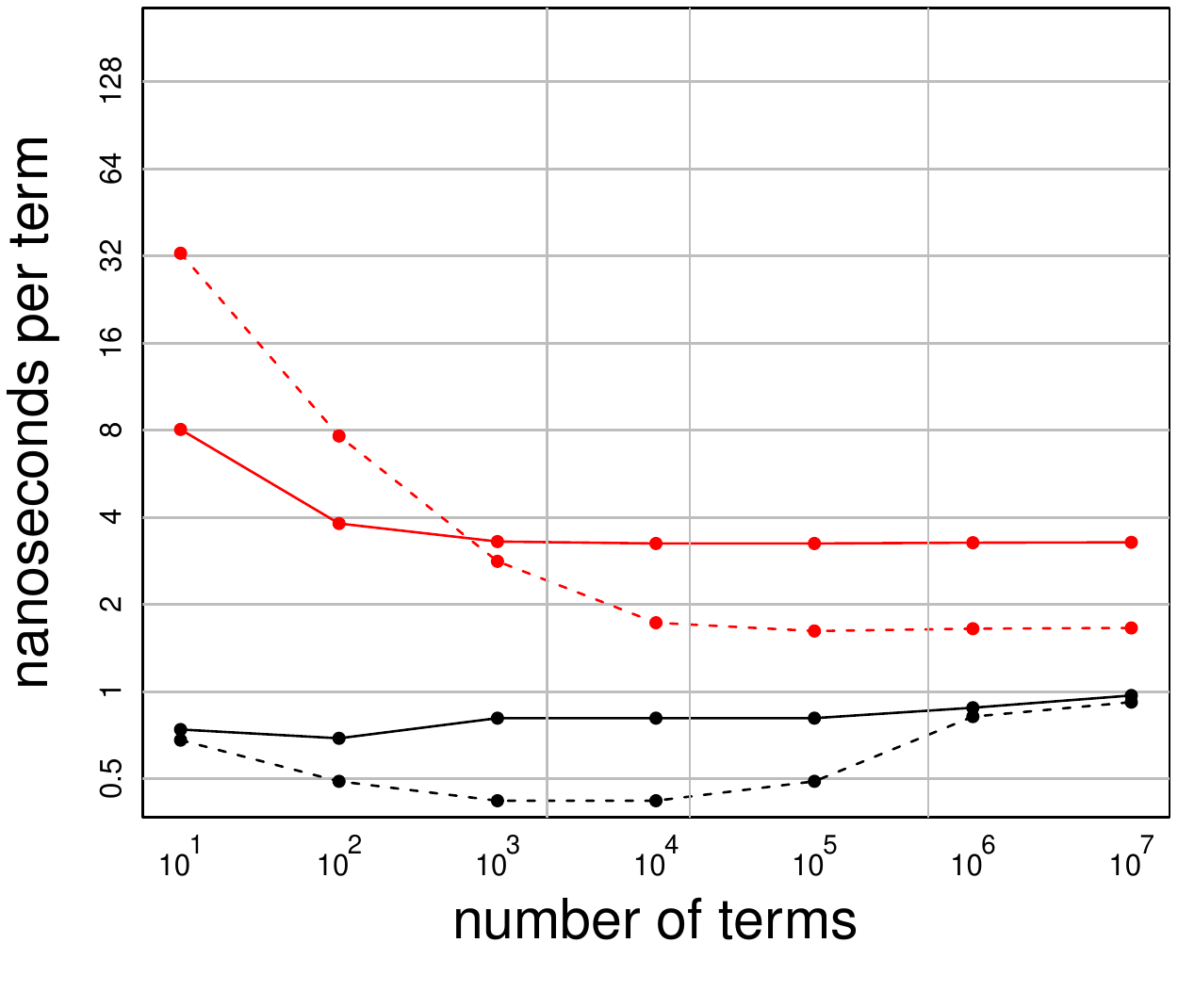} ~~~
   \includegraphics[scale=0.6]{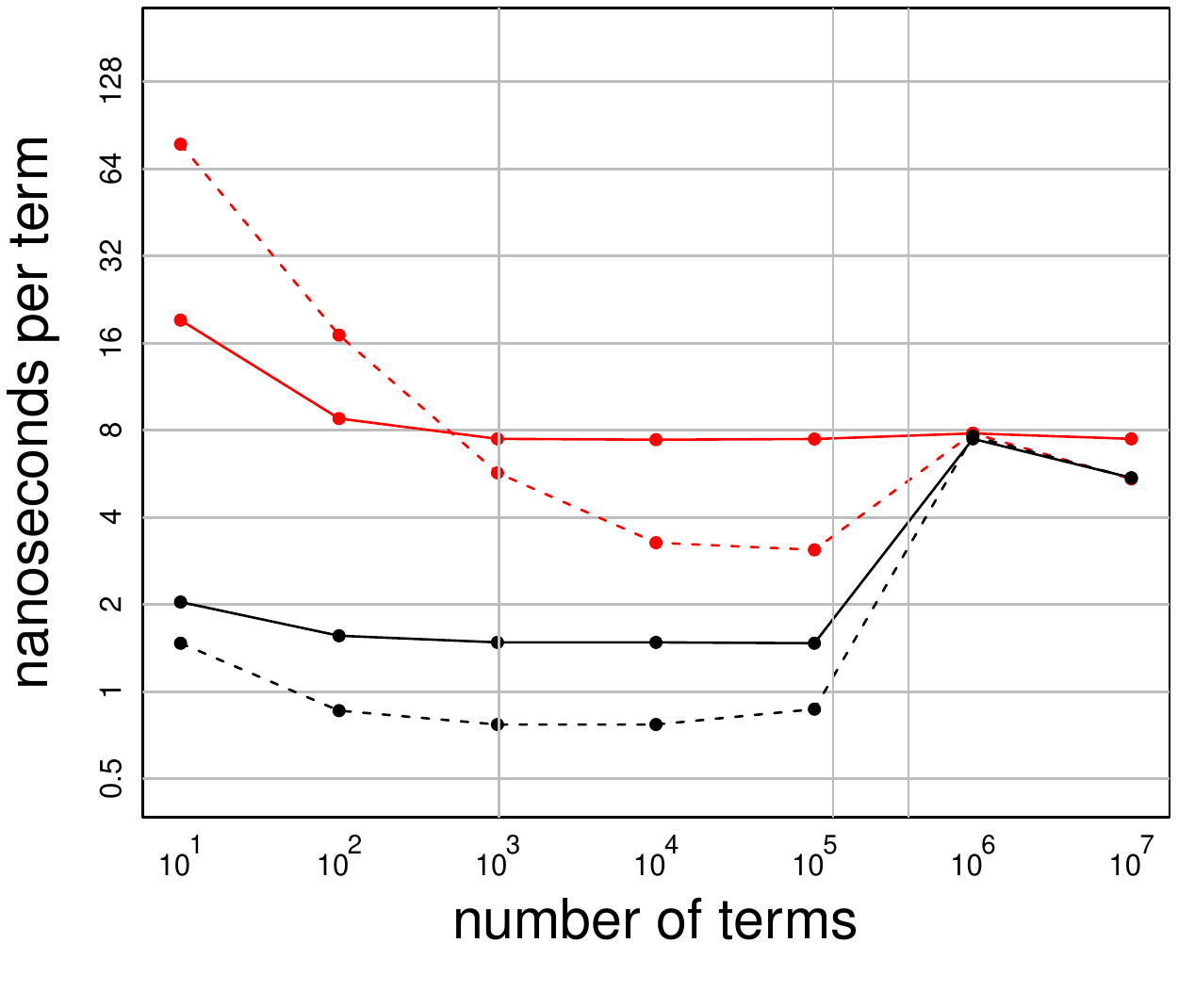}

\vspace*{-11pt}

\caption{Performance of dot product on recent Intel and AMD high-end 
         processors.}\label{fig-res-dot}

\end{figure}

I also implemented functions for computing the squared norm of a
vector (sum of squares of elements) and the dot product of two vectors
(sum of products of corresponding elements) using the small and large
superaccumulator methods for the summations.  (The products were
computed as usual, with rounding to the nearest 64-bit
double-precision floating point number.)  I compared these
implementations with versions using simple ordered and unordered
summation.  The results for the Intel E3-1230~v2 and the AMD Opteron
6348 are shown in Figures~\ref{fig-res-norm} and~\ref{fig-res-dot}.

The times shown in these figures are somewhat disappointing.
Considering that the inner loops of the superaccumulator methods make
no use of the processor's floating-point instructions, I had hoped
that the multiplications in these functions would be executed in
parallel with the integer operations on the superaccumulator, with the
result that the squared norm of a vector would be computed in no more
time than required for summing its elements (and similarly for the dot
product, if the two vectors remain in cache memory).  This is true for
the small superaccumulator method on the AMD Opteron 6348, but for the
large superaccumulator method on this processor, and for both
superaccumulator methods on the Intel E3-1230~v2, the time required
for computing the squared norm is noticeably greater than the time for
summing the elements with the same method.  In contrast the times for
computing the squared norm using simple ordered and unordered
summation are indistinguishable from the times for simple summation,
for vectors of length 1000 or more.  The picture is the same for
computation of the dot product, until the greater memory required
becomes a factor for large vectors.


The reason for this worse than expected performance is not apparent,
but one might speculate that the compilers simply fail to arrange
instructions in a manner that would allow for exploitation of the
instruction-level parallelism that would seem to be possible.  Note,
however, that for large vectors the times to compute the squared norm
or dot product with exact summation are nevertheless still less than a
factor of two greater than the times using simple ordered summation.

More information on these performance assessments, including details
of the computer systems and compilers used, is included in the
supplementary information for this paper.

\subsection*{Discussion}\vspace*{-7pt}

On modern 64-bit processors, serial implementations of the two new
exact summation methods introduced in this paper dominate, in
combination, what appears to be the best combination of previous exact
summation methods --- the iFastSum and OnlineExact methods of Zhu and
Hayes (2010).  The advantage is typically about a factor of two for
large summations.  Note also that the superaccumulator methods produce
a finite final result whenever the correct rounding of the exact sum
is representable as a finite 64-bit floating-point number, whereas the
methods of Zhu and Hayes may produce overflow even when the final
result can be represented.

With the improvement in performance obtained with these
superaccumulator methods, exact summation is less than a factor of two
slower than simple ordered summation, and about a factor of three
slower than simple unordered summation, when summing more than a few
thousand terms.  For large vectors, computing the sum exactly is
faster than attempting to reduce (but not eliminate) error using
Kahan's (1965) method, and Kahan's method has a significant speed
advantage only when the number of terms is less than about one
hundred.

For many applications, the modest extra cost of computing the
exactly-rounded sum may be well worth paying, in return for the
advantages of accuracy.  Exact summation also has the natural
advantage of being reproducible on any computer system that uses
standard floating point, unlike the situation when a variety of
unordered summation methods are used.

The implementation of the small and large superaccumulator methods can
probably be improved.  In the inner loop of the small superaccumulator
method, the conditional branch testing whether a term is positive or
negative could be eliminated (shifting the term right to produce all
0s or all 1s, then XOR'ing to conditionally negate), although this
might be slower when the terms actually all have the same sign.  The
significant variation in performance seen with different compilers may
indicate that none of them are producing close to optimal code.
Future compiler improvements might therefore speed up the performance
of the exact summation methods.  Alternatively, it seems likely that
performance could be improved by rewriting the routines in assembly
language.

One would also expect that using more than one processor core would
allow for faster exact summation.  Collange, Defour, Graillat, and
Iakymchuk (2015a,b) and Chohra, Langlois, and Parello (2015) both
describe parallel implementations of exact summation.  Although these
authors consider a variety of parallel architectures, I will limit
discussion here to parallelizing exact summation on a shared memory
system with multiple general-purpose processor cores or threads.

In this context, any exact summation method can be parallelized in a
straightforward way by simply splitting the array to be summed into
parts, summing each part in parallel (retaining the full exact sum)
and then adding together the partial sums before finally rounding to a
single 64-bit floating point number.  For the methods of this paper,
this would require writing a routine to add together two small
superaccumulators, a straightforward operation that would take time
comparable to that for producing the final rounded result from a small
superaccumulator.  Of course, it is possible that more integrated
algorithms might be somewhat faster, but for large summations, this
simple approach should exploit most of the possible parallelism
available from using a modest number of cores (eg, the two to eight
cores typical on current workstations).

For very large summations, the results in
Figures~\ref{fig-res-intel64} and~\ref{fig-res-amd} suggest that only
a few cores will be needed to reach the limits imposed by memory
bandwidth.  For example, the Intel Xeon \mbox{E3-1230~v2} processor has a
maximum memory bandwidth of 25.6~GBytes/s, which is 0.3125~ns per
8-byte floating-point value.  When summing arrays of $10^7$ elements,
the large superaccumulator method takes 1.16~ns/term, which is 3.7
times larger than the limit imposed by memory bandwidth, suggesting
that 4 cores would be enough to sum terms at the maximum possible
rate.  Since the bandwidth achievable in practice is probably less
than the theoretical maximum, it may be that fewer than 4 cores or
threads would suffice.  (The E3-1230~v2 processor has 4 cores, each of
which can run 2 threads.)

The issue is more complex for summations of around $10^4$ to $10^6$
terms, which may well reside in faster cache memory, which may or may
not be shared between cores.  Experimental evaluations seem essential
to investigating the limits of parallel summation in this regime.

One should note that when comparing methods that all produce the exact
result, and all do so at the maximum rate, limited by memory
bandwidth, the methods can still be distinguished by how many cores or
threads they use in order to achieve this.  This is an important
consideration in the context of a whole application that runs other
threads as well, and in the wider context of a computer system
performing several jobs,

The small superaccumulator method, as well as iFastSum, are rather
slow when summing only a few terms, being ten to twenty times slower
than simple ordered summation.  The small superaccumulator method sets
67 8-byte chunks to zero on initialization, and must scan them all
when producing a rounded result.  This fixed cost dominates the per
term cost when summing only a few terms.  This will limit use of exact
summation in applications where many small summations are done, which
might be of as few as three terms.  (Sums of two terms are exactly
rounded with standard floating-point arithmetic.)

Several approaches could be considered for reducing this fixed
overhead.  One might replace the full array of 67 chunks with a small
list of the non-zero chunks.  Or one might instead keep track of which
chunks are non-zero in a bit array, foregoing actually setting the
value of a chunk until it becomes non-zero, and also using these bits
to quickly locate the non-zero chunks when producing the final rounded
result. These approaches would increase the cost per term, so the
current small superaccumulator method would probably still be the
fastest method for moderate-size summations.

My original motivation for considering exact summation was improving
the accuracy of the sample mean computation in R.  In this
application, the overhead of calling the mean function in the
interpretive R implementation will dominate the fixed overhead of the
small superaccumulator method, so finding a faster method for very
small summations is not essential.

Computing the sample mean by computing the exactly rounded sum of the
data items and then diving by the number of items will not produce the
correct rounding of the exact sample mean, though it will be quite
close (assuming overflow does not occur).  However, it should be
straightforward to write a function that directly produces the correct
rounding of the value in a small superaccumulator divided by a
positive integer.  I plan to soon implement such an exact sample mean
computation in my pqR implementation of R (Neal, 2013--2015).

\subsection*{Acknowledgements}\vspace*{-7pt}

This research was supported by Natural Sciences and Engineering
Research Council of Canada. The author holds a Canada Research Chair
in Statistics and Machine Learning.

\subsection*{References}\vspace*{-7pt}

\leftmargini 0.2in
\labelsep 0in

\begin{description}
\itemsep 2pt

\item
  Chohra, C., Langlois, P., and Parello, D.\ (2015) ``Efficiency of 
  reproducible level 1 BLAS''. \textit{SCAN 2014 Post-Conference Proceedings},
  10 pages.\\
  \verb|http://hal-lirmm.ccsd.cnrs.fr/lirmm-01101723|

\item Collange, S., Defour, D., Graillat, S., and Iakymchuk, R.\ (2015a)
  ``Full-speed deterministic bit-accurate parallel floating-point
  summation on multi- and many-core architectures'', 11 pages.\\
  \verb|https://hal.archives-ouvertes.fr/hal-00949355v2|

\item Collange, S., Defour, D., Graillat, S., and Iakymchuk, R.\ (2015b)
  ``Numerical reproducibility for the parallel reduction on multi- and
   many-core architectures'', 14 pages.\\
   \verb|https://hal.archives-ouvertes.fr/hal-00949355v3|

\item IEEE Computer Society\ (2008) \textit{IEEE Standard for Floating-Point
  Arithmetic}.

\item Kahan, W.\ (1965), ``Further remarks on reducing truncation errors'',
  \textit{Communications of the ACM}, vol.~8, p.~40.

\item Kulisch, U.~W.\ and Miranker, W.~L.\ (1984) ``The arithmetic of
  the digital computer: A new approach'', \textit{SIAM Review}, vol.~28,
  pp.~1--40. 

\item Kulisch, U.\ (2011) ``Very fast and exact accumulation of
  products'', \textit{Computing}, vol.~91, pp.~397-405.

\item Langlois, P., Parello, D., Goossens, B., and Porada, K.\ (2012) ``Less
  hazardous and more scientific research for summation algorithm
  computing times, 34 pages. \\
  \verb|http://hal-lirmm.ccsd.cnrs.fr/lirmm-00737617|

\item Neal, R.~M.\ (2013--2015) pqR --- a pretty quick version of R. 
  \verb|http://pqR-project.org|

\item R Core Team (1995--2015) The R Project for Statistical Computing. 
  \verb|http://r-project.org|

\item Zhu, Y.-K.\ and Hayes, W. B.\ (2010) ``Algorithm 908: Online
  exact summation of floating-point streams'', \textit{ACM
  Transactions on Mathematical Software}, vol.~37, article~37, 13 pages,
  with accompanying software.

\end{description}

\end{document}